\documentclass[amscd,amssymb,verbatim,12pt]{amsart}
\usepackage{amsmath, colonequals}
\usepackage[usenames]{xcolor}
\usepackage[backend=biber, style=alphabetic, maxnames=50, maxalphanames=50]{biblatex}
\usepackage{biblatex}
\usepackage{bm}

\theoremstyle{plain}
\newtheorem{Thm}{Theorem}[section]
\newtheorem{Prop}{Proposition}[section]
\newtheorem{Cor}{Corollary}[section]
\newtheorem{Lem}{Lemma}[section]

\theoremstyle{definition}
\newtheorem{Def}{Definition}[section]

\theoremstyle{remark}
\newtheorem{Rem}{Remark}[section]
\numberwithin{equation}{subsection}

\numberwithin{equation}{section}

\makeatletter
\DeclareFontFamily{OMX}{MnSymbolE}{}
\DeclareSymbolFont{MnLargeSymbols}{OMX}{MnSymbolE}{m}{n}
\SetSymbolFont{MnLargeSymbols}{bold}{OMX}{MnSymbolE}{b}{n}
\DeclareFontShape{OMX}{MnSymbolE}{m}{n}{
    <-6>  MnSymbolE5
   <6-7>  MnSymbolE6
   <7-8>  MnSymbolE7
   <8-9>  MnSymbolE8
   <9-10> MnSymbolE9
  <10-12> MnSymbolE10
  <12->   MnSymbolE12
}{}
\DeclareFontShape{OMX}{MnSymbolE}{b}{n}{
    <-6>  MnSymbolE-Bold5
   <6-7>  MnSymbolE-Bold6
   <7-8>  MnSymbolE-Bold7
   <8-9>  MnSymbolE-Bold8
   <9-10> MnSymbolE-Bold9
  <10-12> MnSymbolE-Bold10
  <12->   MnSymbolE-Bold12
}{}

\let\llangle\@undefined
\let\rrangle\@undefined
\DeclareMathDelimiter{\llangle}{\mathopen}%
                     {MnLargeSymbols}{'164}{MnLargeSymbols}{'164}
\DeclareMathDelimiter{\rrangle}{\mathclose}%
                     {MnLargeSymbols}{'171}{MnLargeSymbols}{'171}
\makeatother

\newcommand {\n}{\noindent}

\newcommand{\C}{{\mathbb C}}
\newcommand{\R}{{\mathbb R}}

\newcommand{\Z}{\mathbb Z}

\newcommand{\LLl}{\mathbb{L}}

\newcommand{\VV}{{\mathcal{V}}}

\newcommand{\Dli}{{\mathcal{D}}'}

\newcommand{\dd}{{\mathrm{d}}}

\newcommand{\vspp}{\vspace*{5pt}}
\newcommand{\lra}{\longrightarrow}

\newcommand{\T}{\mathbb{T}}

\newcommand{\CO}{{\mathcal{O}}}

\newcommand{\dbar}{\overline{\partial}}

\newcommand{\eeps}{\widehat{\otimes}}
\newcommand{\IM}{{\mathsf{Im}\,}}
\newcommand{\RE}{{\mathsf{Re}\,}}

\newcommand{\Ol}{\mathcal{O}}
\DeclareMathOperator{\spn}{span}
\newcommand{\brr}{\overline}
\newcommand{\iprod}{\mathbin{\lrcorner}}
\DeclareMathOperator{\Dif}{Diff}

\DeclareMathOperator{\Aut}{Aut}
\DeclareMathOperator{\ran}{ran}
\begin{document} 

\title[A comparison principle]{A comparison principle between certain Levi-flat compact CR manifolds and systems of real vector fields}
\author{Paulo D. Cordaro \and Vinícius Novelli}

\date{October 9, 2023.
		$^*$Corresponding author: cordaro@ime.usp.br}

\keywords{Levi Flat CR structures, Locally Integrable Structures, Complexes of Differential Operators}

\subjclass[2020]{Primary: 32V10}
\maketitle

\begin{abstract}
We study, in some models of Levi flat CR manifolds, the correspondence between the tangential Cauchy-Riemann complex and the complex defined by the associated real foliation.
\end{abstract}

\section*{Introduction}

 Taking \cite{C} as a starting point, we study a correspondence between certain Levi-flat CR manifolds and the naturally associated real foliation. From the point of view of the theory of locally integrable structures (\cite{T1}, \cite{BCH_book}), one has two different structures defined on the same manifold, and two associated complexes of differential operators: the tangential Cauchy-Riemann complex $\dbar_b$ and the tangential (along the leaves) de Rham complex $\LLl$.  We study the cohomology spaces of these complexes in some models as well as we compare the property of global hypoellipticity in both cases.
\vspp

 In Section 2 we study these questions for Levi flat CR structures defined by locally trivial fiber bundles whose fibers are complex manifolds. In this situation we have at our disposal the Leray spectral sequence which,
 combined with a K\"unneth formula with parameters inspired by \cite{Andreotti1962}, leads to a fair description of the cohomologies for both structures. Among other things we prove that the global solvability for $\dbar_b$ in top forms is equivalent to  the 
 nonexistence of compact components of the fiber (cf. Corollary \ref{23} below).  
\vspp

 A much more involved model is the case of a product $M\times \T^d$ of a compact, complex manifold with the $d$-dimensional torus, where a natural CR Levi-flat structure can be defined after choosing suitable $d$-closed $(0,1)$-forms on $M$. In this case, the associated real structure was studied by several authors (\cite{BCM}, \cite{BCP}, \cite{Araujo2022}, \cite{Araujo2023}).  We introduce such models in Section 3 and describe some of their main properties. We refer the reader to both theorems 
 \ref{main_thm} and \ref{smooth_criterion}, which will be pivotal in what follows.
\vspp

In Section 4, we use a (non-elliptic) laplacian comparison technique to prove equivalences of global hypoellipticity for $\LLl$ and $\dbar_b$. The introduction of the laplacian allows us to define a notion of global hypoellipticity in higher degree forms for these complexes, which coincides with the usual notion on $0$-forms. Our main result in this section is Corollary \ref{73}, which states that if $M$ is balanced in degree $q$ then the global hypoellipticity of $\dbar_b$ and that of $\LLl$
are equivalent in degree $q$.
\vspp

Finally, in Section 5, we study the cohomology of these complexes, again with the laplacian as the main tool. We prove a version of Hodge's theorem (in arbitrary degree) and, in degrees $0$ and $1$, we show an equivalence between global solvability (or validity of the Cousin property for solutions) between $\LLl$ and $\dbar_b$ (cf. Theorem \ref{comp_cohomol_1} below). 
\section{Preliminaries}

\n {\bf (A)} Let $\Omega$ be a smooth, paracompact manifold of dimension $N\geq 1$ which is endowed with a Levi-flat CR structure $\VV$. This means that $\VV$ is smooth involutive subbundle of $ \C T\Omega$ satisfying the following properties:

\begin{enumerate}
\item $\VV \cap\overline \VV=0$;
\item $\VV \oplus \overline \VV$ is an involutive bundle.
 \end{enumerate}

 We shall denote by $n$ the rank of $\VV$ (the value of $n$ is also called the CR dimension of $\VV$) and by $d\doteq N-2n$ the rank of the characteristic set of $\VV$
(recall that the characteristic set of $\VV$ is $\VV^\perp \cap T^*\Omega$, which is real subundle of $T^*\Omega$ since $\VV$ is CR bundle). We assume $d\geq 1$.
\vspp

 Both $\VV$ and $\VV\oplus\overline\VV$ define locally integrable structures on $\Omega$. Indeed for $\VV$ this is result due to L. Nirenberg \cite{Nirenberg} whereas for $\VV\oplus\overline\VV$ it is consequence of Frobenius theorem.
\vspp

\n  {\bf (B)} Any point in $\Omega$ is the center of a coordinate system $(U;x_1,\ldots,x_n,y_1,\ldots,y_n,s_1,\ldots s_d)$ such that $\VV|_U$ is generated by the vectors fields
$\partial/\partial \overline{z}_j$. Notice that then $(\VV\oplus\overline{V})|_U$ is spanned by $\partial/\partial x_j, \partial/\partial y_k$.
\vspp

 The essentially real structure $\VV\oplus \overline{\VV}$ defines a smooth foliation ${\mathcal F}$ in such a way, given any leaf $F\in {\mathcal F}$ (which has dimension $2n$), the restriction of $\VV$ to $F$ defines a complex structure on $F$. In other words each leaf of ${\mathcal F}$ is a complex manifold of dimension $n$.

\n {\bf (C)} Consider the bundles
$$    G_j = \Lambda^j(\C T^*\Omega/\VV^\perp), \quad j=1,\ldots,n; $$
$$   H_j = \Lambda^k(\C T^*\Omega/(\VV\oplus\overline\VV)^\perp), \quad k=1,\ldots,2n. $$
It is well known that the de Rham complex in $\Omega$ induces differential complexes
$$ \overline{\partial}_b : C^\infty(\Omega,G_j)\lra  C^\infty(\Omega,G_{j+1}),\,\, \LLl:  C^\infty(\Omega,H_k)\lra  C^\infty(\Omega,H_{k+1}),$$
whose cohomologies will be denoted respectively by $H^j(\Omega;\overline\partial_b)$ and $H^k(\Omega;\LLl)$.

 Both complexes are locally exact. Hence if we further introduce the sheaf ${\mathcal A}_{\Omega}$ (resp. ${\mathcal B}_\Omega$) of germs of smooth solutions of the equation $\overline{\partial}_b u=$ (resp.  $\LLl u =0$) it follows from standard arguments in sheaf theory \cite{Godement} that
$$    H^j(\Omega,{\mathcal A}_\Omega)=H^j(\Omega;\overline\partial_b), \quad H^k(\Omega,{\mathcal B}_\Omega)=H^k(\Omega;\LLl). $$

 Of special interest is the vanishing of $H^1(\Omega,{\mathcal A}_\Omega)$, for this implies the validity of the first Cousin problem for solutions of $\VV$.

\section{Structures defined by locally trivial fibre bundles}

\n {\bf (A)} In order to study the cohomology of the class of CR structures which we will describe next, it will be necessary to deal with Fréchet sheaves and their completed tensor products. Recall that a sheaf $\mathcal{F}$ over a topological space $X$ is a \textit{Fréchet sheaf} (respectively, \textit{Fréchet-nuclear sheaf}) if the space of sections $\mathcal{F}(U)$ is a Fréchet (respectively, Fréchet-nuclear) space for every open set $U\subset X$ and the restriction maps $\rho_{UV}:\mathcal{F}(U)\to \mathcal{F}(V)$ are continuous for all open sets $V\subset U \subset X$. 

If $\mathcal{F}$ and $\mathcal{G}$ are Fréchet-nuclear sheaves over topological spaces $X$ and $Y$, respectively, we denote by $\mathcal{F}\eeps \mathcal{G}$ the sheaf on $X\times Y$ associated to the presheaf $\mathcal{F}(U)\eeps \mathcal{G}(V)$, $U\subset X$ and $V\subset Y$ open (here, $\eeps$ denotes the completed tensor product in either the $\pi$ or $\varepsilon$ topologies). 

%We shall require a version of the universal coefficient theorem due to \cite{Kaup1967} (see also Grothendieck-Schwartz \cite{GrothendieckSchwartz}). If $E$ is a Fréchet-nuclear space, we denote by $\mathrm{E}$ the constant sheaf on $X$ with stalks $E$.
%\begin{Thm}[\cite{Kaup1967}, Korollar 2]\label{univ_coef1} Let $X$ be a topological space with countable basis, let $\mathcal{F}$ be a Fréchet-nuclear sheaf over $X$ and let $E$ be a Fréchet-nuclear space. Assume that $H^{j}(X,\mathcal{F})$ is Fréchet for all $j=1,\ldots,q$. Then, there is a topological isomorphism
%\begin{equation}\label{univ_coef} 
%H^{q}(X,\mathrm{E}\eeps \mathcal{F}) \simeq E\eeps H^{q}%(X,\mathcal{F}).
%\end{equation}
%\end{Thm}
%\begin{Rem} This result is a consequence of Kaup's version of K\"unneth's formula: let $X$ and $Y$ be topological spaces with countable bases, ${\mathcal F}$ and ${\mathcal G}$ be Fr\'echet-nuclear
%sheaves over $X$ and $Y$ respectively and assume that, for some $q\geq 1$, $H^j(X,{\mathcal F})$ and $H^j(Y,{\mathcal G})$ are Fr\'echet, for $j=1,\dots,q$.
%Then  the space $H^q(X\times Y,\mathcal{F}\eeps\mathcal{G})$ is Fr\'echet and 
%\begin{equation}\label{kunneth_formula}      H^q(X\times %Y,\mathcal{F}\eeps\mathcal{G}) \simeq \bigoplus_{j+k=q} %H^j(X,\mathcal{F})\eeps H^k(Y,\mathcal{G}).\end{equation} 
%ss\end{Rem}

The case which will be of most interest to us is the case where $\mathcal{F}$ is a fine sheaf (namely, the sheaf of smooth functions on a manifold). We record here some basic permanence properties of $\eeps$ when one of the spaces is $C^\infty(X)$, for $X$ a smooth manifold (this space is Fréchet-nuclear for the usual topology, see \cite{TVSDK}).

\begin{Prop}\label{tp_facts} Let $\alpha:F\to F'$ be a continuous linear map between Fréchet-nuclear spaces and let $X$ be a smooth manifold. Let $1\eeps \alpha:C^\infty(X)\eeps F \to C^\infty(X)\eeps F'$ be the induced map on the completed tensor product. Then,
\begin{enumerate}
    \item The kernel of $1\eeps \alpha$ is given by $C^\infty(X)\eeps \ker \alpha$, which is naturally a (closed) subspace of $C^\infty(X)\eeps F$.
    \item If $\alpha$ has closed range, then $1\eeps \alpha$ also has closed range. More precisely,
    \[
    \ran (1\eeps \alpha) = C^\infty(X) \eeps \ran \alpha.
    \]
\end{enumerate}
    
\end{Prop}
\begin{proof} The second item follows from the fact that the functor $E\mapsto C^\infty(X)\eeps E$ preserves short exact sequences of Fréchet-nuclear spaces (this remains true replacing $C^\infty(X)$ by any Fréchet-nuclear space). See, for instance, p. 435 in \cite{Demailly} or p. 205 in \cite{Andreotti1962} for a direct proof in the case of $C^\infty(X)$.

This property does not (directly) imply item $(1)$, because there we do not assume $\alpha$ has closed range. Here, we take advantage of the fact that $C^\infty(X)$ allows for the use of Fourier series (this argument is based on Andreotti-Grauert \cite{Andreotti1962}): indeed, if $X=\T^d$ is the $d$-dimensional torus, then the result is immediate from Fourier series decomposition:
\[
\ker (1\eeps \alpha) = \left\{\sum_{n\in \Z^d}f_n e^{in\theta} \in C^\infty(\T^d)\eeps F;\,\alpha(f_n)=0\text{ for all }n \right\} \simeq C^\infty(\T^d) \eeps \ker \alpha.
\]
If $X$ is a general (paracompact) smooth manifold, let $\{U_i\}$ be a locally finite covering of $X$ by coordinate charts $h_i:U_i \to \R^d$ such that $h_i(U_i)$ is an open subset of the unit cube $I_d \subset \R^d$. Let $\{\rho_i\}$ be a smooth partition of unity subordinated to this covering. If $f\in C^\infty(X)\eeps F$ is such that $(1\eeps \alpha)(f) = 0$, then $f_i:=\rho_i f$ is also in the kernel of $1\eeps \alpha$. It is clear that $f_i$ can be identified with an element in $C^\infty(\T^d)\eeps F$, which yields $f=\sum f_i \in C^\infty(X)\eeps \ker \alpha$. The reverse inclusion is clear, since the algebraic tensor product $C^\infty(X)\otimes \ker \alpha$ is clearly contained in the (closed) kernel of $1\eeps \alpha$.
\end{proof}

\n {\bf (B)}  We shall consider locally trivial fibre bundles  $(\Omega,\Omega\stackrel{f}\lra S,M)$, where
$\Omega$, the total space, is a smooth manifold of dimension $N$; $S$, the base space, is a smooth manifold of dimension $d$ and $M$, the fibre space,
is a complex manifold of complex dimension $n$. We take the structure group of this bundle to be $\Aut(M)$ (the group of biholomorphisms of $M$). It is easily seem that a natural Levi flat CR structure, of CR dimension equal to $n$, can be introduced on $\Omega$. Indeed, if $U\subset S$ is an open set  over which $(\Omega,\Omega\stackrel{f}\lra S,M)$ trivializes, that is,
there is a smooth diffeomorphism $h: U\times M\simeq f^{-1}(U)$ satisfying
$(f\circ h)(x,z)=x$, $x\in U$ and $z\in M$, it follows that $h_\star(\{0\}\times T^{1,0}M)$ defines a Levi flat CR structure of CR dimension equal to $n$
on $f^{-1}(U)$ and, when we vary $U$, they all match together to define
a Levi flat CR structure of CR dimension equal to $n$ on $\Omega$.
\vspp

We fix a (countable) cover $\mathcal{U}:=(U_j)_{j\in \Z_+}$ of $S$ such that the bundle is trivialized in each $U_j$. We denote it by $\phi_j:U_j \times M \to f^{-1}(U_j)$, with cocycle map $\phi_{ij}:=\phi_i^{-1}\circ \phi_j: (U_i \cap U_j) \times M \to (U_i \cap U_j)\times M$ given by $\phi_{ij}(x,p)=(x,g_{ij}(x)\cdot p)$, where $g_{ij}:U_i \cap U_j \to \Aut(M)$ is smooth, for $i,j\in \Z_+$ such that $U_i \cap U_j \not=\emptyset$. The following lemma computes the cohomology of $\mathcal{A}_{\Omega}$ over such a set.

\begin{Lem}\label{lem1_ag} Let $U_i$ be a trivializing open set. Then, if $H^q(M,\Ol_M)$ is Hausdorff for some $q\geq 1$, $H^{q}(f^{-1}(U_i),\mathcal{A}_{\Omega}\big|_{f^{-1}(U_i)})$ is also Hausdorff and there is a topological isomorphism
\[
T_i:H^{q}(f^{-1}(U_i),\mathcal{A}_{\Omega}\big|_{f^{-1}(U_i)}) \rightarrow C^\infty(U_i) \eeps H^q(M,\Ol_M).
\]
Moreover, this isomorphism is compatible with restrictions to open subsets of $U_i$ (that is, it induces a sheaf isomorphism).
\end{Lem}
\begin{proof} If $\phi_i:U_i \times M \to f^{-1}(U_i)$ is the trivialization, then the pullback sheaf $\phi_i^{\ast} \mathcal{A}_{\Omega}\big|_{f^{-1}(U_i)}$ is the sheaf $C^{\infty}_{U_i}\eeps \Ol_M$. In particular, the cohomology is isomorphic to the cohomology of the complex
\[
C^\infty(U_i\times M;\Lambda^{q-1}) \xrightarrow{\dbar_{q-1}} C^\infty(U_i\times M;\Lambda^q) \xrightarrow{\dbar_q}C^\infty(U_i\times M;\Lambda^{q+1}).
\]
of $q$-forms on $M$ depending smoothly on parameters in $U_i$. This is equivalent to the complex
\begin{equation}\label{local_complex}
C^\infty(U_i)\eeps C^\infty(M;\Lambda^{q-1}) \xrightarrow{1\eeps \dbar_{q-1}} C^\infty(U_i)\eeps C^\infty(M;\Lambda^{q}) \xrightarrow{1\eeps \dbar_q} C^\infty(U_i)\eeps C^\infty(M;\Lambda^{q+1}).
\end{equation}
We conclude from Proposition \ref{tp_facts} that the cohomology of \ref{local_complex} is isomorphic to
\[
\frac{C^\infty(U_i)\eeps \ker \dbar_{q}}{C^\infty(U_i) \eeps \ran \dbar_{q-1}} \simeq C^\infty(U_i) \eeps H^q(M,\Ol_M),
\]
since the (completed) tensor product is an exact functor and $H^q(M,\Ol_M)$ is a Fréchet space.
\end{proof}

Assume that $H^{q}(M,\Ol_M)$ is Hausdorff for some $1\leq q \leq n$. Then, we consider the sequence space 
\[
X:= \prod_{j\in \Z_+}C^\infty\left(U_j, H^{q}(M,\Ol_M)\right) \simeq \prod_{j\in \Z_+}C^{\infty}(U_j) \eeps H^{q}(M,\Ol_M),
\]
with its natural structure of Fréchet space. We define a natural subspace of $X$ by
\begin{equation}\label{trivializations_1}
\mathcal{E}(\mathcal{U},M,\Ol_M):=\left\{(v_i)_{i\in \Z_+}\in X; \,\left(T_i \circ T_j^{-1}\right) v_j\big|_{U_i \cap U_j} = v_i\big|_{U_i \cap U_j}  \right\}.
\end{equation}
It is clear that this is a closed subspace of $X$. With these notations in mind, we can state the main result.

\begin{Thm}\label{app_prop1} Let $(\Omega,\Omega\stackrel{f}\longrightarrow S,M)$ be as above. Let $\mathcal{U}=(U_j)_{j\in \Z_+}$ be a trivializing cover of $S$ and assume that $H^{q}(M,\Ol_M)$ is Hausdorff for some $1\leq q \leq n$. Then, there is a topological isomorphism
\begin{equation}\label{cohomology_isom1}
H^{q}(\Omega,\mathcal{A}_{\Omega}) \simeq \mathcal{E}(\mathcal{U},M,\Ol_M).
\end{equation}\end{Thm}
\vspp

\begin{proof} The key ingredient in the proof is the use of the Leray spectral sequence associated to $(\Omega,\Omega\stackrel{f}\longrightarrow S,M)$.
Let $R^q f_{\ast}\mathcal{A}$ denote the $q$-th cohomology sheaf of the direct image sheaf $f_{\ast}\mathcal{A}$, that is, the sheaf associated to the presheaf
$$ V\mapsto H^q(f^{-1}(V),{\mathcal A}_\Omega), \mbox{\quad $V\subset S$,\, open.}$$
 It is well known that there exists a spectral sequence (the Leray spectral sequence)
$$  E^{p,q}_2 = H^p(S,R^q f_{\ast}\mathcal{A}) \Longrightarrow H^{p+q}(\Omega,\mathcal{A}_\Omega).$$

%% In order to compute the term $E^{p,q}_2$ we shall use a special covering of $S$. Indeed, thanks to (\cite{BottTu}, pp. 42-43) we can construct an open
%covering ${\mathcal U}=\{U_i\}_{i\in I}$ of $S$ with the following properties:
%\begin{enumerate}
%\item[(i)] $(\Omega,\Omega\stackrel{f}\longrightarrow S,M)$ trivializes over $U_i$, for every $i\in I$;
%\item[(ii)] Given $q\geq 1$ and distinct indices $i_0,\ldots i_q\in I$
%the intersection $U_{i_0}\cap\cdots\cap U_{i_q}$ is diffeomorphic to $\R^d$.
%\end{enumerate}

 Our first observation is that $R^q f_{\ast}\mathcal{A}$ and
$C^\infty_S\eeps \mathbf{H}^q(M,\CO_M)$ coincide over each $U_i$ (here we denote by $\mathbf{H}^q(M,\CO_M)$ the constant sheaf on $S$ with stalks
$H^q(M,\CO_M)$): indeed, since sheafification is a functor (\cite{KashiwaraSchapira}, p. 85), this follows from Lemma \ref{lem1_ag}.

\vspp

In particular, we conclude that $R^{q}f_{\ast}\mathcal{A}_{\Omega}$ is locally fine, and therefore, locally soft. Then, it has no higher cohomology (see, for instance, Proposition 4.13 and Theorem 4.15, pages 204 and 205 in \cite{Demailly}). We obtain that 
\[
E^{p,q}_2 = 0,\,\,p\geq 1,
\]
and, therefore, the Leray spectral sequence degenerates. This implies a topological isomorphism
\[
H^{q}(\Omega,\mathcal{A}_{\Omega}) \simeq H^{0}(S,R^{q} f_{\ast}\mathcal{A}_{\Omega}).
\]
Observe that we have an embedding 
\begin{align*}
\Phi: H^{0}(S,R^{q} f_{\ast}\mathcal{A}_{\Omega}) &\hookrightarrow \prod_{j\in \Z_+}R^q f_{\ast}\mathcal{A}_\Omega(U_i)\\
\sigma &\mapsto \left(\sigma\big|_{U_i} \right)_{i\in \Z_+}
\end{align*}
whose range is given by
\[
\IM \Phi = \left\{(u_i)_{i\in \Z_+}\in\prod_{j\in \Z_+}R^q f_{\ast}\mathcal{A}_\Omega(U_i) ;\,\,u_i\big|_{U_i \cap U_j} = u_j\big|_{U_i\cap U_j} \right\}.
\]
Applying the isomorphisms from Lemma \ref{lem1_ag}, we obtain an isomorphism of Fréchet spaces
\[
\Psi:H^{q}(\Omega,\mathcal{A}_{\Omega}) \to \mathcal{E}(\mathcal{U},M,\Ol_M).
\]
\end{proof}

 By a similar argument, now recalling the fact that the spaces $H^q(M,\C)$ are always Fr\'echet, thanks to de Rham theorem, we can obtain an analogous description of the space $H^{q}(\Omega,\mathcal{B}_{\Omega})$ (compare with Theorem 14.18 in \cite{BottTu}). Note that in this case, the structure group can be the full group $\Dif(M)$ of diffeomorphisms of $M$, and the complex structure of the fibre is not required. Indeed, introducing the closed subspace
 \[
 \mathcal{E}(\mathcal{U},M,\C) \subset Y:=\prod_{j\in \Z_{+}}
C^\infty(U_j)\eeps H^q(M,\C)
\]
in exactly the same way as \ref{trivializations_1}, but now considering the trivializations induced by the real structure on the bundle, we can state the following
\begin{Thm}\label{app_prop2} Let $(\Omega,\Omega\stackrel{f}\longrightarrow S,M)$ be as above. Let $\mathcal{U}=(U_j)_{j\in \Z_+}$ be a trivializing cover of $S$. Then, there is a topological isomorphism
\begin{equation}\label{cohomology_isom2}
H^{q}(\Omega,\mathcal{B}_{\Omega}) \simeq \mathcal{E}(\mathcal{U},M,\C)
\end{equation}
for all $1\leq q \leq 2n$.
\end{Thm}

 We  now list some consequences of these results.

 \begin{Cor}\label{P}
  Let $1\leq q\leq n$. Then $H^q(M,\mathcal{O}_M)=0$ if and only if 
 $H^q(\Omega,\mathcal{A}_\Omega)=0$. In particular, if $M$ is Stein then $H^q(\Omega,\mathcal{A}_\Omega)=0$ for $q=1,\dots,n$ and, from Theorem
 \ref{app_prop2}, $H^k(\Omega,\mathcal{B}_\Omega)=0$ for $k=n+1,\dots,2n$.
 \end{Cor}

 \begin{proof} By Theorem \ref{app_prop1} $H^q(M,\mathcal{O}_M)=0$ implies
 $H^q(\Omega,\mathcal{A}_\Omega)=0$. Conversely, assume that $H^q(\Omega,\mathcal{A}_\Omega)=0$ and let  $\omega\in C^\infty_{(0,q)}(M)$
 be $\overline{\partial}$-closed. Let $U\subset S$ be an open set such that there is a diffeomorphism $h: U\times M\simeq f^{-1}(U)$ satisfying $(f\circ h)(x,z)=x$, $x\in U$ and $z\in M$. Select $\psi\in C^\infty_c(U)$ which is equal to one at some point $x_0\in U$ . Then $(h^{-1})^*(\psi\otimes \omega)$, extended as zero outside $f^{-1}(U)$, defines an element $\beta \in C^\infty(\Omega,G_q)$ which is $\overline{\partial}_b$-closed. Since
 $H^q(\Omega,\mathcal{A}_\Omega)=0$ there exists $\alpha\in
C^\infty (\Omega, G_{q-1})$ such that $\overline{\partial}_b \alpha= \beta$.
Then $\alpha_\bullet \doteq h^*(\alpha|_{f^{-1}U})$ solves $\overline{\partial}\alpha_\bullet = \psi\otimes \omega$ in $U\times M$.
In particular $\gamma \doteq\alpha_\bullet (x_0,\cdot)\in C^\infty_{(0,q-1)}(M)$ solves $\overline{\partial} \gamma = \omega$ in $M$.

The remaining statements follow from well known results for Stein manifolds.
\end{proof}

\vspp
\begin{Cor}\label{23} The cohomology space $H^n(\Omega,{\mathcal A}_\Omega)$ is trivial if and only if $M$ has no compact connected component.

\end{Cor}
\begin{proof} By Corollary \ref{P}, $H^n(\Omega,{\mathcal A}_\Omega)=0$ if and only if $H^n(M,\CO_M)=0$, and this last space is trivial if and only if
$M$ has no compact component (cf. \cite{Malgrange}).

\end{proof}

\begin{Cor} Assume the bundle $\Omega \to S$ is trivial. Then, if $H^{q}(M,\Ol_M)$ is Hausdorff for some $1\leq q \leq n$, we have a topological isomorphism
\[
H^{q}(\Omega,\mathcal{A}_{\Omega}) \simeq C^\infty(S) \eeps H^{q}(M,\Ol_M).
\]
\end{Cor}
\begin{Rem} Observe that this is a generalization of Proposition 7, page 208 in \cite{Andreotti1962}. Moreover, a similar result holds for $H^{q}(\Omega,\mathcal{B}_{\Omega})$ (see the observations that precede Theorem \ref{app_prop2}).
    
\end{Rem}

 %By a similar argument, now recalling the fact that the spaces $H^q(M,\C)$ are always Fr\'echet, thanks to de Rham theorem, we obtain a generalization of Theorem 14.18 in \cite{BottTu}.  Note that in this case, the structure group can be the full group $\Dif(M)$ of diffeomorphisms of $M$, and the complex structure of the fibre is not required.

%\begin{Thm}\label{app_prop2} Let $(\Omega,\Omega\stackrel{f}\longrightarrow S,M)$ be as above. If $1\leq q\leq 2n$ then
%\[ H^q(\Omega;{\mathcal B}_\Omega) \simeq C^\infty(S)\eeps H^q(M,\C). \]

%\end{Thm}

% We  now list some consequences of Theorems \ref{app_prop1} and \ref{app_prop2}:

%\begin{Cor} Each of the following hypothesis implies \eqref{cohomology_isom1}:

%\begin{itemize}
%\item $M$ is compact;
%\item $M$ is a Stein manifold (in this case $H^q(\Omega, {\mathcal %A}_\Omega)=0$ for $q=1,\ldots,n$ and then, by Theorem %\ref{app_prop2}, $H^k(\Omega, {\mathcal B}_\Omega)=0$ for
%$k=n+1,\ldots, 2n$).
%\end{itemize}

%\end{Cor}
%\vspp
%\begin{Cor}\label{23} The cohomology space $H^n(\Omega,{\mathcal A}_\Omega)$ is trivial if and only if $M$ has no compact connected component.

%\end{Cor}
%\begin{proof} Indeed, it is well known that $H^n(M,\Ol_M)$ is always a Fr\'echet space and that it is trivial if and only if $M$ has no compact connected component. 

%\end{proof}

%\begin{Cor}  If $H^1(M,\Ol_M)=0$ then $H^1(\Omega,{\mathcal A}_\Omega)=0$.
%
%\end{Cor}

%\begin{Cor} Suppose that $M$ is K\"ahler and compact. Then
%\[ H^q(\Omega,{\mathcal A}_\Omega)\simeq H^q(\Omega,{\mathcal %B}_\Omega), \quad q=0,\ldots,n. \]

%\end{Cor}

\section{A class of compact Levi-flat CR structures}

\n {\bf (A)} In this part of the work, we concentrate on a different class of Levi-flat CR structures in which the complex structure of the (leaves of the) foliation is not fixed, but varies in a controlled way. We start by establishing some notation. Let $M$ be a compact, connected complex manifold of (complex) dimension $n\geq 1$ and let $\T^d = \R^d/\Z^d$ be the $d$-dimensional torus. Consider $\omega_1,\ldots,\omega_d \in C^\infty(M;\Lambda^{0,1}T^{\ast}M)$ forms of type $(0,1)$ such that $\dbar \omega_k = 0$ for all $k=1,\ldots,d$. Let $(\theta_1,\ldots,\theta_k)$ denote the usual angular coordinates on the torus. We define 
\[
T':=T^{1,0}M \oplus \spn\{\alpha_1,\ldots,\alpha_d\},
\]
where $\alpha_k:=\dd \theta_k + \omega_k$ for $k=1,\ldots,d$ (we are identifying forms and bundles on $M$ and on $\T^d$ with their pullbacks to $M\times \T^d$).
\begin{Prop} $T'$ defines a locally integrable CR structure on $M\times \T^d$. 
\end{Prop}
\begin{proof} $T'$ is clearly a subbundle of $\C T^{\ast}(M\times \T^d)$, which is locally integrable from Dolbeault-Grothendieck's lemma. Moreover, we have $T'+\brr{T'}=\C T^{\ast}M$, since $\dd \theta_k=\alpha_k-\omega_k$ is a section of $T'+\brr{T'}$ for every $k=1,\ldots,d$, which shows $T'$ is a CR structure.
\end{proof}
\n {\bf(B)} Denote by $\mathcal{V}\subset \C T(M\times \T^d)$ the subbundle orthogonal to $T'$ (for the duality between one-forms and vector fields), which has rank $n$. Consider a system of holomorphic coordinates $(z_1,\ldots,z_n)$ in an open set $\Omega \subset M$, and write $\omega_k = \sum_{j=1}^n\omega_{jk}\dd \brr{z_j}$, for $k=1,\ldots,d$, where $\omega_{jk} \in C^\infty(\Omega)$. In such coordinates, a frame for the the bundle  $\mathcal{V}$ over $\Omega \times \T^d$ is given by the set of vector fields
\[
L_j = \frac{\partial}{\partial \brr{z_j}} - \sum_{k=1}^d \omega_{jk} \frac{\partial}{\partial \theta_k},\,\,j=1,\ldots,n.
\]
We shall compute the characteristic set and the Levi form of $\mathcal{V}$. Fixing a point $p=(z,\theta) \in \Omega \times \T^d$, let $v\in T'_p$. Then, we can write
\[
v=\sum_{j=1}^n \alpha_j \dd z_j + \sum_{k=1}^d \xi_k \alpha_k = \sum_{j=1}^n \alpha_j \dd z_j + \sum_{k=1}^d \xi_k \dd \theta_k + \sum_{j=1}^n \sum_{k=1}^d \xi_k \omega_{jk}\dd \brr{z_j},\,\,\alpha_j,\xi_k \in \C.
\]
This covector is real if and only if $\xi_k \in \R$ for all $k=1,\ldots,d$ and 
\[
\alpha_j = \sum_{k=1}^d \xi_k \brr{\omega_{jk}},\,\,j=1,\ldots,n.
\]
We conclude that the characteristic set $T^{\circ}=T'\cap T^\ast(M\times \T^d)$ has (real) dimension $d$, and is generated by the real forms $\dd \theta_k + 2\RE \omega_k$, $k=1,\ldots,d$. Computing the Lie brackets $[L_j,\brr{L_k}]$, we obtain
\[
[L_j,\brr{L_k}] = \sum_{l=1}^d \left(\frac{\partial \omega_{jl}}{\partial z_k} - \brr{\frac{\partial \omega_{kl}}{\partial z_j}} \right)\frac{\partial}{\partial \theta_l},\,\,j,k=1,\ldots,n.
\]
Therefore, the matrix of the Levi form $\mathfrak{L}_{(p,v_k)}$ with respect to the basis $\{L_1,\ldots,L_n\}$ at the characteristic vector $v_l=\dd \theta_l + 2\RE \omega_l$ is given by
\[
\mathfrak{L}_{(p,v_l)} = \frac{1}{2i}\left(\frac{\partial \omega_{jl}}{\partial z_k} - \brr{\frac{\partial \omega_{kl}}{\partial z_j}} \right)_{1\leq j,k \leq n}.
\]
We obtain the following result:
\begin{Prop}\label{prop1} Let $\omega_1,\ldots,\omega_d$ be smooth $(0,1)$-forms on $M$ satisfying $\dbar \omega_k = 0$ for all $k=1,\ldots,d$. Then, the following are equivalent:
\begin{enumerate}
\item $T'$ is Levi-flat.
\item $\dd(\omega_k+\brr{\omega_k})=0$ for all $k=1,\ldots,d$.
\end{enumerate}

\end{Prop}
\begin{proof} Just observe that, in a system of holomorphic coordinates $(z_1,\ldots,z_N)$ where $\omega_l = \sum_{j=1}^N \omega_{jl}\dd \brr{z_j}$, we have
\[
\dd(\omega_l+\brr{\omega_l}) = \sum_{j,k=1}^n \left(\frac{\partial \omega_{jl}}{\partial z_k} - \brr{\frac{\partial \omega_{kl}}{\partial z_j}} \right)\dd \brr{z_k}\wedge \dd z_j,
\]
since $\dbar \omega_l = 0$.
\end{proof}

\n {\bf (C)} From now on, we assume the structure $T'$ is Levi-flat ($\dd(\omega_k+\brr{\omega}_k)=0$ for all $k=1,\ldots,d$). The tangential CR complex associated to such a structure the following: for $0\leq q \leq n$,
\begin{equation}\label{cr_complex}
\left(\dbar_b\right)_q:C^\infty\left(M\times \T^d;\Lambda^{q}\C T^{\ast}M\right) \to C^\infty\left(M\times \T^d;\Lambda^{q+1}\C T^\ast M\right)
\end{equation}
given by
\[
\dbar_b u = \dbar u - \sum_{k=1}^d \omega_k \wedge  \frac{\partial u}{\partial \theta_k},
\]
where the derivatives $\partial/\partial \theta_k$ are defined component-wise. One can also consider this complex acting on distributional sections, i.e., currents on $M\times \T^d$ (valued in $\C T^{\ast}M$). We would like to deduce regularity properties of solutions of $\dbar_b u = f$ from the corresponding properties of solutions of the real complex
\[
\LLl_q:C^\infty\left(M\times \T^d; \Lambda^{q}\C T^{\ast}M\right) \to C^\infty\left(M\times \T^d;\Lambda^{q+1}\C T^{\ast}M\right)
\]
given by
\begin{equation}\label{real_complex}
\LLl_q u = \dd_{M} u - \sum_{k=1}^d (\omega_k + \brr{\omega_k}) \wedge \frac{\partial u}{\partial \theta_k}.
\end{equation}
(this defines a differential complex since $\dd(\omega_k+\brr{\omega_k})=0$ for all $k=1,\ldots,d$). The main technique we will use to compare both complexes involves the Laplacians, and we discuss this next.

\n {\bf {(D)}} First, we discuss the complex acting on functions (i.e., $q=0$). Fix an hermitian metric on $M$ (with corresponding volume form denoted by $\dd V$) and consider the usual flat metric on the torus $\T ^d$. This metric induces hermitian products in the exterior algebra $\Lambda^\bullet(\C T^{\ast}M)$. We can then consider the Hilbert spaces $L^2(M;\dd V)$ (respectively, $L^2(M\times \T^d;\dd V \dd \theta)$) and $L^2(M,\C T^\ast M;\dd V)$ (respectively, $L^2(M\times \T^d;\C T^\ast M)$) given by the (equivalence classes of) measurable sections of the indicated bundles that satisfy
\[
\int_M \langle f,f\rangle \dd V<\infty,\,\,(\text{respectively, }\int_{M\times \T^d}\langle f,f\rangle \dd V \dd \theta < \infty).
\]
We shall use the notation $\llangle \cdot,\cdot \rrangle$ for the global hermitian product and $\langle \cdot,\cdot \rangle(x)$ for the product in $\Lambda^\bullet(\C T^\ast_pM)$. We can then consider the Hilbertian adjoints
\[
\dbar_b^\ast,\LLl^\ast:C^\infty(M\times \T^d;\C T^\ast M)\to C^\infty(M\times \T^d).
\]
We shall compute these adjoints explicitely:
\begin{Prop}\label{prop2} Let $\omega_1,\ldots,\omega_d \in C^\infty(M\times \T^d;\Lambda^{0,1}M)$ be $(0,1)$-forms satisfying $\dbar \omega_k = 0$ for all $k=1,\ldots,d$. Then,
\begin{enumerate}
\item The adjoint of the map $\omega_k \frac{\partial}{\partial \theta_k}:C^\infty(M\times \T^d)\to C^\infty(M\times \T^d;\C T^\ast M)$ is the map $\beta \mapsto -\left \langle \frac{\partial \beta}{\partial \theta_k},\omega_k \right \rangle$.
\item $\langle \dd_{\,M} u,\omega_k+\brr{\omega_k}\rangle(x)=\langle \partial u,\brr{\omega_k}\rangle(x)+\langle \dbar u,\omega\rangle(x)$ for all $u\in C^\infty(M\times \T^d)$, $x\in M\times \T^d$ and $k=1,\ldots,d$ (the same formula holds for the global product $\llangle \cdot,\cdot \rrangle$).
\item $|\omega_k+\brr{\omega_k}|^2(p)=2|\omega_k|^2(p)$ for all $p\in M$ and $k=1,\ldots,d$ (the same formula holds for the global product $\llangle \cdot,\cdot \rrangle$).
\end{enumerate}

\end{Prop}
\begin{proof} For item $1)$, let $f\in C^\infty(M\times \T^d)$ and $\beta \in C^\infty(M\times \T^d;\C T^\ast M)$. Then,
\begin{align*}
\left \llangle \omega_k \frac{\partial f}{\partial \theta_k},\beta \right \rrangle &= \int_{M\times \T^d}\left \langle \frac{\partial f}{\partial \theta_k}(z,\theta)\omega_k(z),\beta(z)\right \rangle\dd V\dd \theta \\
&=\int_{M\times \T^d}\frac{\partial f}{\partial \theta_k}(z,\theta)\langle \omega_k(z),\beta(z)\rangle \dd V \dd \theta \\
&=-\int_{M\times \T^d}f(z,\theta)\left \langle \omega_k(z), \frac{\partial \beta}{\partial \theta_k}\right \rangle \dd V \dd \theta \\
&=\left \llangle f, -\left\langle\frac{\partial \beta}{\partial \theta_k},\omega_k \right\rangle \right \rrangle.
\end{align*}
For the items $2)$ and $3)$, just observe that forms with different bidegrees are orthogonal and that $\langle \brr{v},\brr{w}\rangle(x) =\brr{\langle v,w \rangle(x)}$.
\end{proof}
We shall now perform some computations in local coordinates. Let $\Omega \subset M$ be an open subset with holomorphic coordinates $(z_1,\ldots,z_n)$. We write the hermitian metric in these coordinates as 
\[
h=\sum_{j,k=1}^n h_{jk}\dd z_j \otimes \dd \brr{z_k},
\]
where $(h_{jk})$ is a hermitian matrix of functions in $C^\infty(\Omega)$ (which are given by $h_{jk}=\langle \partial/\partial z_j,\partial/\partial z_{k}\rangle$). We denote the inverse of this matrix by $(h^{jk})$ (it's a simple exercise to verify that $h^{jk}=\langle \dd z_k,\dd z_j\rangle$). Let $\phi,\psi \in C^\infty(\Omega;\Lambda^{0,1}M)$, which we write as
\[
\phi = \sum_{j=1}^N \phi_j \dd \brr{z_j},\,\,\psi = \sum_{k=1}^N \psi_k \dd \brr{z_k},\,\,\phi_j,\psi_k \in C^\infty(\Omega).
\]
Then, we have
\begin{align*}
\langle \phi,\psi\rangle(z)&= \sum_{j,k=1}^N \phi_j(z) \brr{\psi_k}(z) \langle \dd \brr{z_j},\dd \brr{z_k}\rangle \\
&=\sum_{j,k=1}^n \phi_j(z) \brr{\psi_k}(z) h^{jk}(z).
\end{align*}
In the same way, if $\phi$ e $\psi$ are of type $(1,0)$ (given by $\phi=\sum \phi_j \dd z_j$ and $\psi = \sum \psi_j \dd z_j$), we get
\[
\langle \phi, \psi \rangle(z) = \sum_{j,k=1}^{n}\phi_j(z)\brr{\psi_k}(z)h^{kj}(z)=\sum_{j,k=1}^{n}\phi_j(z)\brr{\psi_k}(z)\brr{h^{jk}(z)}.
\]
We write, for fixed $k=1,\ldots,d$, 
\[
\omega_k = \sum_{j=1}^N \omega_{jk}\dd \brr{z_j},\,\,\omega_{jk}\in C^\infty(\Omega).
\]
We would like to compare $d^{\ast}(f(\omega_k+\brr{\omega_k}))$ with $\dbar^{\ast}(f\omega_k)$ (where $f\in C^\infty(M\times \T^d)$), under the hypothesis that $\omega_k + \brr{\omega_k}$ is closed. Let $u\in C^\infty_c(\Omega\times \T^d)$. Then,
\begin{align*}
\llangle u,\dd^{\ast}\left(f(\omega_k+\brr{\omega_k}) \right)-2\dbar^{\ast}(f\omega_k)\rrangle &= \llangle \dd u,f\omega_k+ f\brr{\omega_k} \rrangle - 2\llangle \dbar u,f\omega_k \rrangle \\
&=\llangle \partial u,f\brr{\omega_k}\rrangle - \llangle \dbar u,f\omega_k\rrangle.
\end{align*}
Writing $h(z)=\det (h_{jk}(z))$ gives us (see, for instance, page 146 in \cite{Kodaira})
\begin{align*}
&\llangle u,\dd^{\ast}(f(\omega_k+\brr{\omega_k})) - 2\dbar^{\ast}(f\omega_k)\rrangle=\\
&= 2^n\int_{U\times \T^d}\left\{\sum_{j,l=1}^n\frac{\partial u}{\partial z_j}(z,\theta)\brr{f(z,\theta)}\omega_{lk}(z)\brr{h^{jl}(z)}-\frac{\partial u}{\partial \brr{z_j}}(z,\theta)\brr{f(z,\theta)\omega_{lk}(z)}h^{jl}(z) \right\}h(z)\dd x \dd \theta \\
&= 2^n \sum_{j,l=1}^n\int_{U\times \T^d}\left\{\frac{\partial u}{\partial z_l}(z,\theta) \omega_{jk}(z) - \frac{\partial u}{\partial \brr{z_j}}(z,\theta)\brr{\omega_{lk}(z)} \right\}h^{jl}(z)h(z)\brr{f(z,\theta)}\dd x \dd \theta.
\end{align*}
Now we shall integrate by parts. When we integrate the terms of the form $\omega_{lk}$, the term that will appear will be of the form $\sum_{j,l}\frac{\partial \omega_{jk}}{\partial z_l}-\brr{\frac{\partial \omega_{lk}}{\partial z_j}}$, which vanishes since $\dd(\omega_k + \brr{\omega_k})=0$ (see the proof of Proposition \ref{prop1}). Therefore, we need only to integrate the terms of the form $h^{jl}h$ (which we call (I)) and $\brr{f}$ (which we call (II)). We obtain
\begin{align*}
\text{(I)}&=-2^n\sum_{j,l=1}^n \int_{U\times \T^d}\left\{\frac{\partial (h^{jl}h)}{\partial z_l}\omega_{jk}(z)-\frac{\partial (h^{jl}h)}{\partial \brr{z_j}}\brr{\omega_{lk}(z)} \right\}u(z,\theta)\brr{f(z,\theta)}\dd x \dd \theta \\
&=2^n  \int_{U\times \T^d}u(z,\theta) \left\{\brr{f(z,\theta)} \frac{1}{h(z)}\sum_{j,l=1}^n \left(\frac{\partial (h^{jl}h)}{\partial \brr{z_j}}\brr{\omega_{lk}(z)} - \frac{\partial (h^{jl}h)}{\partial z_l}\omega_{jk}(z) \right) \right\} h(z)\dd x \dd \theta \\
&=\int_{U\times \T^d}u(z,\theta)\cdot \brr{\left\{f(z,\theta) \cdot \left(\frac{1}{h(z)}\sum_{j,l=1}^n \frac{\partial (\brr{h^{jl}h})}{\partial z_j}\omega_{lk}(z) - \frac{\partial (\brr{h^{jl}h})}{\partial \brr{z_l}}\brr{\omega_{jk}(z)} \right) \right\}} 2^n h(z)\dd x \dd \theta \\
&=\llangle u,f(R_k-\brr{R_k})\rrangle,
\end{align*}
where $R_k\in C^\infty(\Omega)$ is given by
\[
R_k(z)=\frac{1}{h(z)}\sum_{j,l=1}^n \frac{\partial (\brr{h^{jl}h})}{\partial z_j}\omega_{lk}(z).
\]
For the term (II), we have
\begin{align*}
\text{(II)}&=-2^n\sum_{j,l=1}^n \int_{U\times \T^d}\left\{\frac{\partial \brr{f}}{\partial z_l}(z,\theta)\omega_{jk}(z)- \frac{\partial \brr{f}}{\partial \brr{z_j}}(z,\theta)\brr{\omega_{lk}(z)} \right\}h^{jl}(z)h(z)u(z,\theta)\dd x\dd \theta \\
&=2^{n}\int_{U\times \T^d}u(z,\theta)\brr{\left\{\sum_{j,l=1}^{n}\frac{\partial f}{\partial z_j}(z,\theta)\omega_{lk}(z)\brr{h^{jl}(z)}-\sum_{j,l=1}^n\frac{\partial f}{\partial \brr{z_l}}(z,\theta)\brr{\omega_{jk}(z)} \brr{h^{jl}(z)}\right\}}h(z)\dd x \dd \theta \\
&=\llangle u,\langle \partial f,\brr{\omega_k} \rangle - \langle \dbar f,\omega_k\rangle \rrangle.
\end{align*}
Since $u\in C^\infty_c(\Omega \times \T^d)$ is arbitrary, we conclude that 
\begin{equation}\label{eq1}
\dd^{\ast}\left(f(\omega_k+\brr{\omega_k}) \right)-2\dbar^\ast(f\omega_k) = \langle \partial f,\brr{\omega_k} \rangle - \langle \dbar f,\omega_k \rangle + f(R_k-\brr{R_k})
\end{equation}
in $\Omega \times \T^d$. With these results in hand, we compare the Laplacians $\dbar_b^{\ast}\dbar_b$ and $\LLl^\ast \LLl$ (when acting on functions supported on $\Omega \times \T^d$). From Proposition \ref{prop2}, we have 
\[
\dbar_b^{\ast} = \dbar^{\ast} + \sum_{k=1}^d \left \langle \frac{\partial}{\partial \theta_k},\omega_k \right \rangle,\,\,\LLl^{\ast}=\dd^{\ast}+\sum_{k=1}^d \left \langle \frac{\partial}{\partial \theta_k},\omega_k+\brr{\omega_k}\right \rangle.
\]
We obtain then
\begin{align*}
\dbar_b^\ast \dbar_b &= \left( \dbar^{\ast} + \sum_{k=1}^d \left \langle \frac{\partial\,\cdot}{\partial \theta_k},\omega_k \right \rangle\right)\left(\dbar-\sum_{k=1}^d \omega_k \frac{\partial}{\partial \theta_k} \right) \\
&=\dbar^\ast \dbar - \sum_{k=1}^d \dbar^{\ast}\left(\omega_k \frac{\partial}{\partial \theta_k} \right) + \sum_{k=1}^d \left \langle \frac{\partial}{\partial \theta_k}\dbar,\omega_k \right \rangle - \sum_{k,k'=1}^d \langle \omega_{k'},\omega_{k}\rangle \frac{\partial^2}{\partial \theta_k \partial \theta_{k'}}.
\end{align*}
In the same way,
\begin{align*}
\LLl^{\ast}\LLl &= \left(\dd^{\ast}+\sum_{k=1}^d \left \langle \frac{\partial}{\partial \theta_k},\omega_k+\brr{\omega_k}\right \rangle \right)\left(\dd - \sum_{k=1}^d (\omega_k+\brr{\omega_k})\frac{\partial}{\partial \theta_k} \right) \\
&=\dd^{\ast}\dd - \sum_{k=1}^d \dd^{\ast}\left((\omega_k+\brr{\omega_k})\frac{\partial}{\partial \theta_k} \right)+\sum_{k=1}^d \left \langle \frac{\partial}{\partial \theta_k} \dd, \omega_k + \brr{\omega_k}\right \rangle - \sum_{k,k'=1}^d \langle \omega_{k'}+\brr{\omega_{k'}},\omega_k+\brr{\omega_k}\rangle \frac{\partial^2}{\partial \theta_k \partial \theta_{k'}} \\
&=\dd^{\ast}\dd - \sum_{k=1}^d \dd^{\ast}\left((\omega_k+\brr{\omega_k})\frac{\partial}{\partial \theta_k} \right)+\sum_{k=1}^d \left( \left \langle\frac{\partial}{\partial \theta_k}\partial,\brr{\omega_k}\right \rangle + \left \langle\frac{\partial}{\partial \theta_k}\dbar,\omega_k \right \rangle\right) - 2\sum_{k,k'=1}^d \langle \omega_{k'},\omega_{k}\rangle \frac{\partial^2}{\partial \theta_k \partial \theta_{k'}}
\end{align*}
Putting everything together, we obtain
\begin{align*}
\LLl^\ast \LLl - 2\dbar^\ast \dbar &=\dd^{\ast}\dd - 2\dbar^{\ast}\dbar - \sum_{k=1}^d \left(\dd^{\ast}\left((\omega_k+\brr{\omega_k})\frac{\partial}{\partial \theta_k} \right)-2\dbar^{\ast}\left(\omega_k \frac{\partial}{\partial \theta_k} \right) \right) +\\
&+ \sum_{k=1}^d \left(\left \langle \frac{\partial}{\partial \theta_k}\partial,\brr{\omega_k} \right \rangle  - \left \langle\frac{\partial}{\partial \theta_k}\dbar,\omega_k \right \rangle\right) \\
&=\dd^{\ast}\dd - 2\dbar^{\ast}\dbar - \sum_{k=1}^d \left(\left \langle \partial \frac{\partial}{\partial \theta_k},\brr{\omega_k} \right \rangle - \left \langle \dbar \frac{\partial}{\partial \theta_k},\omega_k \right \rangle + (R_k-\brr{R_k})\frac{\partial}{\partial \theta_k}\right) + \\
&+\sum_{k=1}^d \left(\left \langle \frac{\partial}{\partial \theta_k}\partial,\brr{\omega_k} \right \rangle  - \left \langle\frac{\partial}{\partial \theta_k}\dbar,\omega_k \right \rangle\right) \\
&=\dd^{\ast}\dd - 2 \dbar^{\ast}\dbar - \sum_{k=1}^d (R_k-\brr{R_k})\frac{\partial}{\partial \theta_k},
\end{align*}
in $\Omega \times \T^d$. Taking $f\equiv 1$ in \eqref{eq1} shows that $R_k-\brr{R_k}=\dd^{\ast}(\omega_k+\brr{\omega_k}) - 2\dbar^\ast(\omega_k) = \partial^{\ast}\brr{\omega_k}-\dbar^\ast \omega_k$. Therefore, we can invariantly write the identity
\begin{equation}\label{eq2}
\LLl^{\ast}\LLl - 2\dbar^{\ast}_b\dbar_b = \Box_{\dd}-2\Box_{\dbar} + 2i\sum_{k=1}^d(\IM \dbar^{\ast}\omega_k)\frac{\partial}{\partial \theta_k},
\end{equation}
where $\Box_{\dd}=\dd^{\ast}\dd$ and $\Box_{\dbar}=\dbar^{\ast}\dbar$ are the Laplace-Beltrami and Laplace-Dolbeault's operators, acting on functions (this follows because both these operators are local and are equal in small neighborhoods of every point of $M\times \T^d$).

If the manifold $M$ admits a particular kind of metric, this expression can be simplified even further.
\begin{Def}\label{def1} Let $0\leq q \leq 2n$. A complex manifold $M$ is \textit{balanced in degree} $q$ if it admits a hermitian metric $h$ such that $\frac{1}{2}\Box_{\dd}^q:=\frac{1}{2}\left(\dd^{\ast}_{\,q}\dd_{\,q} + \dd_{\,q-1} \dd^{\ast}_{\,q-1} \right)=(\dbar^{\ast}_q\dbar_q + \dbar_{q-1} \dbar^{\ast}_{q-1})=:\Box_{\dbar}^q$. 
\end{Def}
\n It is well-known that K\"ahler manifolds are balanced in every degree $0\leq q \leq 2n$ (see Corollary 6.5, page 306 in \cite{Demailly}). Moreover, a result of \cite{Hsiung} shows that, conversely, if a complex manifold is balanced in degree $0$ and $1$, then it is a K\"ahler manifold. However, for fixed values of $q$, the class of manifolds that are balanced in that particular degree might include non-K\"ahler manifolds (for example, in dimension $n\geq 3$, there are balanced manifolds in degree zero that are not K\"ahler \cite{Michelsohn}).
\begin{Thm}\label{teo6} Assume the manifold $M$ is balanced in degree $0$. Then, 
\[
\frac{1}{2}\LLl^{\ast}\LLl =\dbar^\ast_b \dbar_b,
\]
when the adjoints are taken with respect to the balanced metric.
\end{Thm}
\begin{proof} Let $k=1,\ldots,d$ be fixed. Since $\dd(\omega_k+\brr{\omega_k}) = 0$, we can write in a sufficiently small neighborhood $\Omega$ of a point $z_0 \in M$ the equation $\omega_k + \brr{\omega_k}=\dd f$, where $f\in C^\infty(\Omega;\R)$, by Poincaré's lemma. In particular, we have $\dbar f = \omega_k$. Then, in $\Omega$,
\[
\dbar^{\ast}\omega_k = \dbar^\ast \dbar f = \frac{1}{2}\Box_{\dd}f,
\]
which is a real-valued function (again, using $\Box_{\dd}=2\Box_{\dbar}$). The identity now follows from \ref{eq2}.
\end{proof}

\n {\bf(E)} Now we shall study what happens in larger degree (i.e., $q\geq 1$). We consider now the Hilbert spaces $L ^2(M;\Lambda^{p,q}M;\dd V);$ and $L^2(M\times \T^d;\Lambda^{p,q};\dd V \dd \theta)$, for $0\leq p,q \leq n$. We shall recall some elementary geometric operations.
\begin{Def} Let $N$ be a smooth manifold and $X \in \C T_pN$ be a tangent vector at $p\in N$. If $u\in \Lambda^{q}T^{\ast}_pN$ is a $q$-covector over on $p$ ($q\geq 1$), then the contraction $X\iprod u$ is defined by 
\[
(X\iprod u)(v_1,\ldots,v_{q-1}):=u(X,v_1,\ldots,v_{q-1}),\,\,v_j \in \C T_p N.
\]
\end{Def}
\begin{Def} Let $(N,g)$ be an $n$-dimensional smooth Riemannian manifold. The musical isomorphism at $p\in N$ is defined as
\begin{align*}
\flat_p:\C T_p N &\rightarrow \C T_p^{\ast}N \\
v&\mapsto \left(\C T_p N \ni w \mapsto \langle w,v \rangle_p \right),
\end{align*}
with inverse denoted by $\sharp_p:=\flat_p^{-1}:\C T^{\ast}_pN\to \C T_p N$.
\end{Def}
\begin{Rem} Let $(x_1,\ldots,x_n)$ be a system of coordinates in $N$. Then, writing $g_{ij}=\langle \partial/\partial x_i,\partial/\partial x_j\rangle$ for the coefficients of the Riemannian metric, we have the following: let $X=\sum_{i=1}^n X_i \frac{\partial}{\partial x_i}$ be a tangent vector. Then, if $X^{\flat} = \sum_{i=1}^n \alpha_i \dd x_i$, 
\begin{align*}
\alpha_i = X^{\flat}(\partial/\partial x_i) &= \left \langle \frac{\partial}{\partial x_i},X \right \rangle \\
&=\sum_{j=1}^{n}g_{ij}\brr{X_j},
\end{align*}
i.e., $X^{\flat}=\sum_{i=1}^n \left(\sum_{j=1}^n g_{ij}\brr{X_j} \right) \dd x_i$. In the same way, if $\omega=\sum_{i=1}^n \omega_i \dd x_i$, we have $\omega^{\sharp} = \sum_{i=1}^n \left(\sum_{j=1}^n g^{ij}\brr{\omega_j} \right)\frac{\partial}{\partial x_i}$, where $(g^{ij})$ is the inverse matrix of $(g_{ij})$.
\end{Rem}
\n Now, if $\omega \in C^\infty(N;\C T^{\ast}N)$ is a one-form in $N$, we can define the interior product by $\omega$ using the musical isomorphism: if $\alpha \in C^\infty(N;\Lambda^qT^{\ast}N)$, $q\geq 1$, then
\[
\omega \iprod \alpha := \omega^{\sharp}\iprod \alpha.
\]
\begin{Lem}\label{lemma1} Let $(N,g)$ be a Riemannian manifold and let $p\in N$. 
\begin{enumerate}
\item Let $\omega_p \in CT_p^{\ast}N$, $\alpha_p \in \Lambda^{q}T^{\ast}_pN$ and $\eta_p \in \Lambda^{q-1}T^{\ast}_p N$, $q\geq 1$. Then,
\[
\langle \omega_p \iprod \alpha_p,\eta_p \rangle= \langle \alpha_p,\omega_p \wedge \eta_p \rangle.
\]
In particular, if $\omega \in C^\infty(N;\C T^{\ast}N)$, $\alpha \in C^\infty(N;\Lambda^qT^{\ast}N)$ and $\eta \in C^\infty(N;\Lambda^{q-1}T^{\ast}N)$, we have $\llangle \omega \iprod \alpha,\eta \rrangle = \llangle \alpha,\omega \wedge \eta \rrangle$.
\item The derivation property
\[
\omega_p \iprod (u_p \wedge v_p) = (\omega_p \iprod u_p)\wedge v_p + (-1)^{q} u_p \wedge (\omega_p \iprod v_p)
\]
holds, for $\omega_p \in \C T^{\ast}_p N$, $u_p \in \Lambda^{q}T^{\ast}_p N$ and $v_p \in \Lambda^{r}T^{\ast}_p N$, $q+r \geq 1$. 
\item If $\omega_1,\omega_2 \in \C T^{\ast}_p N$, then
\[
\omega_1 \iprod \omega_2 = \langle \omega_2,\omega_1\rangle.
\]
\end{enumerate}
\end{Lem}
\begin{proof} Item $1):$ Let $\{\xi_1,\ldots,\xi_n\}$ be an orthonormal basis for $\C T_p N$. Then, $\{\xi_1^{\flat},\ldots,\xi_n^{\flat}\}$ is the dual basis, and is also an orthonormal basis for $\C T^{\ast}_p N$. The result is clear for $\omega_p = \xi_i^{\flat}$, $\alpha_p = \xi^{\flat}_I$ and $\eta_p = \xi^{\flat}_J$, $|I|=q$, $|J|=q-1$, and this implies the general case by linearity. Item $2)$ is easily verified using this orthonormal basis. For the final item, writing $\omega_i = \sum_{j=1}^{n}\omega_{ij}\xi_j^{\flat}$, $i=1,2$, we have 
\[
\omega_1 \iprod \omega_2 = \sum_{l,k=1}^{n}\brr{\omega_{1l}} \omega_{2k}\xi_l \iprod \xi_k^{\flat} = \sum_{l=1}^n \brr{\omega_{1l}}\omega_{2l}=\langle \omega_2,\omega_1 \rangle.
\]
\end{proof}
\begin{Rem} We remark that the interior product is local, i.e., if $U\subset M$ is an open set, $(\omega \iprod \alpha)\big|_{U}=(\omega\big|_{U})\iprod (\alpha\big|_{U})$.

\end{Rem}
\n We prove a basic identity for $\dd^{\ast}$ and $\dbar^{\ast}$.
\begin{Lem}[Leibniz formula]\label{leibniz} Let $U\subset M$ be an open set, $\alpha \in C^\infty(U;\Lambda^{q}T^{\ast}M)$ be a $q$-form (with $q \geq 1$) and $f\in C^\infty(U)$. Then,
\begin{equation}\label{leibniz_d}
\dd^{\ast}\left(f\alpha \right) = f \left(\dd^{\ast}\alpha\right) - \left(\dd \brr{f}\right)\iprod \alpha\text{ in }U
\end{equation}
and
\begin{equation}\label{leibniz_dbar}
\dbar^{\ast}\left(f\alpha \right) = f\left(\dbar^{\ast} \alpha\right) - \left(\dbar \brr{f}\right)\iprod \alpha\text{ in }U.
\end{equation}
\end{Lem}
\begin{proof} We prove only \eqref{leibniz_d} (the other one being analogous). Let $\beta \in C^{\infty}_c(U;\Lambda^{q-1}T^{\ast}M)$. Then,
\begin{align*}
\llangle \dd^{\ast}(f\alpha),\beta\rrangle &= \llangle \alpha, \brr{f}\dd \beta\rrangle \\
&=\llangle \alpha, \dd(\brr{f}\beta) - (\dd \brr{f})\wedge \beta \rrangle \\
&=\llangle f(\dd^{\ast}\alpha) - (\dd \brr{f})\iprod \alpha,\beta \rrangle.
\end{align*}
\end{proof}
\n Now we can proceed to the main result.
\begin{Thm}\label{main_thm} Assume that $M$ is balanced in degree $q \in \{0,\ldots,2N\}$. Then,
\begin{equation}\label{main_eq}
\frac{1}{2}\Box_{\LLl}^{q}= \Box_{\dbar_b}^q.
\end{equation}
\end{Thm}
\begin{proof} The previous section covers the case $q=0$, so we assume $q\geq 1$. Fix a point $p\in M$ and let $U_p \subset M$ be an open neighborhood of $p$ such that there exist real-valued functions $f_1,\ldots,f_d \in C^\infty(U_p;\R)$ such that $\dd f_k =\omega_k+\brr{\omega_k}$ for all $k=1,\ldots,d$ (this neighborhood exists by Poincaré's lemma, since $\dd(\omega_k+\brr{\omega_k})=0$ for all $k=1,\ldots,d$). We shall prove that 
\[
\frac{1}{2}\Box^q_{\LLl} u = \Box_{\dbar_b}^{q}u
\]
for all $u\in C^\infty_c(U_p \times \T^d;\Lambda^{q}T^{\ast}M)$. Note that, by applying a partition of unity, this implies that $\frac{1}{2}\Box_{\LLl}^q = \Box^q_{\dbar_b}$. Let $u \in C^\infty_c(U_p\times \T^d;\Lambda^{q}\C T^{\ast}M)$. Then,
\begin{align*}
\Box_{\LLl}^q u &= \LLl^{\ast}_q \LLl_q u + \LLl_{q-1}\LLl^{\ast}_{q-1}u \\
&=\LLl^{\ast}_q \left(\dd u - \sum_{k=1}^{d}(\omega_k+\brr{\omega_k})\wedge \frac{\partial u}{\partial \theta_k} \right) + \LLl_{q-1}\left(\dd^{\ast}u + \sum_{k=1}^{d}(\omega_k+\brr{\omega_k}) \iprod \frac{\partial u}{\partial \theta_k}\right) \\
&=\dd^{\ast}\dd u - \sum_{k=1}^{d}\dd^{\ast}\left((\omega_k+\brr{\omega_k})\wedge \frac{\partial u}{\partial \theta_k} \right) + \sum_{k=1}^{d}(\omega_k+\brr{\omega_k})\iprod \frac{\partial (\dd u)}{\partial \theta_k} \\
&-\sum_{k,k'=1}^{d}(\omega_k+\brr{\omega_k}) \iprod \left((\omega_{k'}+\brr{\omega_{k'}})\wedge \frac{\partial^2 u}{\partial \theta_k \partial \theta_{k'}} \right) +\dd \dd^{\ast}u +\sum_{k=1}^{d}\dd \left((\omega_k+\brr{\omega_k})\iprod \frac{\partial u}{\partial \theta_k} \right) \\
&-\sum_{k=1}^{d}(\omega_k+\brr{\omega_k})\wedge \frac{\partial (\dd^{\ast} u)}{\partial \theta_k} - \sum_{k,k'=1}^d (\omega_k+\brr{\omega_k}) \wedge \left((\omega_{k'}+\brr{\omega_{k'}})\iprod \frac{\partial^2 u}{\partial \theta_k \partial \theta_{k'}} \right).
\end{align*}
From Lemma \ref{lemma1},
\begin{align*}
&\sum_{k,k'=1}^{d} (\omega_k+\brr{\omega_k}) \iprod \left((\omega_{k'}+\brr{\omega_{k'}})\wedge \frac{\partial^2 u}{\partial \theta_k \partial \theta_{k'}}\right) =\\
& =\sum_{k,k'=1}^d \left(\left(\omega_k+\brr{\omega_k} \right)\iprod \left(\omega_{k'}+\brr{\omega_{k'}} \right) \right) \wedge \frac{\partial^2 u}{\partial \theta_k \partial \theta_{k'}} - \left(\omega_{k'}+\brr{\omega_{k'}} \right) \wedge \left((\omega_k+\brr{\omega_k})\iprod \frac{\partial^2 u}{\partial \theta_k \partial \theta_{k'}} \right),
\end{align*}
which implies that (again, using \ref{lemma1} and the orthogonality relations)
\begin{align*}
&\sum_{k,k'=1}^{d} (\omega_k+\brr{\omega_k}) \iprod \left((\omega_{k'}+\brr{\omega_{k'}})\wedge \frac{\partial^2 u}{\partial \theta_k \partial \theta_{k'}}\right)+\sum_{k,k'=1}^d (\omega_k+\brr{\omega_k}) \wedge \left((\omega_{k'}+\brr{\omega_{k'}})\iprod \frac{\partial^2 u}{\partial \theta_k \partial \theta_{k'}} \right) = \\
&=2\sum_{k,k'=1}^d  \left(\omega_k\iprod \omega_{k'}\right) \wedge \frac{\partial^2 u}{\partial \theta_k \partial \theta_{k'}}.
\end{align*}
Now we exploit the locality of the operators. From \eqref{leibniz}, we have
\begin{align*}
\sum_{k=1}^d \Box_{\dd}^{q} \left(f_k \frac{\partial u}{\partial \theta_k} \right) &= \sum_{k=1}^d \left(\dd \left(f_k \frac{\partial \dd^{\ast}u}{\partial \theta_k} - (\omega_k+\brr{\omega_k})\iprod \frac{\partial u}{\partial \theta_k} \right)+\dd^{\ast}\left((\omega_k+\brr{\omega_k})\wedge \frac{\partial u}{\partial \theta_k} + f_k\frac{\partial \dd u}{\partial \theta_k} \right)\right) \\
&=\sum_{k=1}^d \left((\omega_k+\brr{\omega_k}) \wedge \frac{\partial \dd^{\ast}u}{\partial \theta_k} + f_k \frac{\partial \dd \dd^{\ast}u}{\partial \theta_k}  \right) - \sum_{k=1}^d \dd \left((\omega_k+\brr{\omega_k})\iprod \frac{\partial u}{\partial \theta_k} \right) +\\
&+\sum_{k=1}^d \dd^{\ast}\left((\omega_k+\brr{\omega_k})\wedge \frac{\partial u}{\partial \theta_k} \right) + \sum_{k=1}^d f_k \frac{\partial \dd^{\ast}\dd u}{\partial \theta_k} - \sum_{k=1}^d (\omega_k+\brr{\omega_k})\iprod \frac{\partial \dd u}{\partial \theta_k}.
\end{align*}
We obtain, then,
\begin{align*}
\sum_{k=1}^d \left(\Box_{\dd}^{q} \left(f_k \frac{\partial u}{\partial \theta_k} \right) - f_k \frac{\partial (\Box_d^{q}u)}{\partial \theta_k}\right) &=\sum_{k=1}^d (\omega_k+\brr{\omega_k})\wedge \frac{\partial \dd^{\ast}u}{\partial \theta_k} - \sum_{k=1}^d \dd \left((\omega_k+\brr{\omega_k})\iprod \frac{\partial u}{\partial \theta_k} \right) \\
&+\sum_{k=1}^d \dd^{\ast}\left((\omega_k+\brr{\omega_k})\wedge \frac{\partial u}{\partial \theta_k} \right) -\sum_{k=1}^d (\omega_k+\brr{\omega_k})\iprod \frac{\partial \dd u}{\partial \theta_k}.
\end{align*}
Going back to the expression of $\Box_{\LLl}^q u$, we obtain the fundamental identity
\begin{align*}
\Box_{\LLl}^q u = \Box_{\dd}^q u + \sum_{k=1}^{d}\left(f_k \frac{\partial(\Box_{\dd}^q u)}{\partial \theta_k} - \Box_{\dd}^q \left(f_k\frac{\partial u}{\partial \theta_k} \right)\right) - 2\sum_{k,k'=1}^{d}(\omega_k\iprod \omega_{k'})\wedge \frac{\partial^2 u}{\partial \theta_k \partial \theta_{k'}}.
\end{align*}
In the same way, for $\Box_{\dbar_b}^q = (\dbar_b)_q^{\ast}(\dbar_b)_q + (\dbar_b)_{q-1}(\dbar_b)^{\ast}_{q-1}$ we have
\begin{align*}
\Box_{\dbar_b}^qu &= (\dbar_b)_q^{\ast}(\dbar_b)_qu + (\dbar_b)_{q-1}(\dbar_b)^{\ast}_{q-1}u \\
&=(\dbar_b)_q^{\ast}\left(\dbar u - \sum_{k=1}^d \omega_k \wedge \frac{\partial u}{\partial \theta_k} \right) + (\dbar_b)_{q-1}\left( \dbar^{\ast}u+\sum_{k=1}^d \omega_k \iprod \frac{\partial u}{\partial \theta_k}\right) \\
&=\dbar^{\ast}\dbar u -\sum_{k=1}^{d}\dbar^{\ast}\left(\omega_k \wedge \frac{\partial u}{\partial \theta_k} \right) + \sum_{k=1}^d \omega_k \iprod \frac{\partial (\dbar u)}{\partial \theta_k} - \sum_{k,k'=1}^d \omega_{k} \iprod \left(\omega_{k'}\wedge \frac{\partial^2 u}{\partial \theta_k \partial \theta_{k'}} \right) + \\
&+\dbar \dbar^{\ast}u + \sum_{k=1}^d \dbar\left(\omega_k \iprod \frac{\partial u}{\partial \theta_k} \right) - \sum_{k=1}^d \omega_k \wedge \frac{\partial (\dbar^{\ast}u)}{\partial \theta_k} - \sum_{k,k'=1}^d \omega_k \wedge \left(\omega_{k'}\iprod \frac{\partial^2 u}{\partial \theta_k \partial \theta_{k'}} \right).
\end{align*}
Making exactly the same maneuvers we did for $\Box_{\LLl}^q$ yields the identity (recall that $\dbar f_k = \omega_k$ for all $k=1,\ldots,d$)
\[
\Box_{\dbar_b}^q u = \Box_{\dbar}^q u + \sum_{k=1}^{d}\left(f_k \frac{\partial(\Box_{\dbar}^q u)}{\partial \theta_k} - \Box_{\dbar}^q \left(f_k\frac{\partial u}{\partial \theta_k} \right)\right) - \sum_{k,k'=1}^{d}(\omega_k\iprod \omega_{k'})\wedge \frac{\partial^2 u}{\partial \theta_k \partial \theta_{k'}}.
\]
The hypothesis that $\frac{1}{2}\Box_{\dd}^q = \Box_{\dbar}^q$ yields the result.
\end{proof}
\n {\bf (F)} We will, in this section, determine criteria that allow us to decide whether or not a differential form $f$ on $M\times \T^d$ (valued in $\C T^\ast M$) is smooth by analyzing the decay (in $L^\infty$ or $L^2$ norms) of the Fourier coefficients $\widehat{f}(j)$, $ j\in \Z^d$. To do so, we need to consider the operators $\LLl$ and $\dbar_b$ acting on currents: we use the notation $\Dli(M\times \T^d;\Lambda^{p}\C T^\ast M)$ for the space of currents on $M\times \T^d$ of degree $p$, valued in $\C T^\ast M$. Any such current can be written locally a differential $p$-form on $M$ with distributional coefficients defined on $M\times \T^d$. 

We define the following spaces, for $1\leq q \leq 2n$:
\[
X_P^{q} =\left\{\alpha \in \Dli(M\times \T^d;\Lambda^{q}\C T^\ast M);\,P_q\alpha \text{ and }P^\ast_{q-1} \alpha \text{ are smooth sections} \right\}.
\]
where $P$ is either $\LLl$ or $\dbar_b$. If $\alpha \in X_P^q$, then we can write a Fourier decomposition
\[
\alpha = \sum_{j\in \Z^d}\widehat{\alpha}(j)e^{i\langle j,\theta \rangle},\,\theta \in \T^d,
\]
where $\widehat{\alpha}(j) \in \Dli(M;\Lambda^q \C T^{\ast}M)$ is a $q$-current in $M$.  A simple verification shows that
\[
P \alpha = \sum_{j\in \Z^d}P_j (\widehat{\alpha}(j)) e^{i\langle j,\theta \rangle},
\]
where $P_j$ is $D_j:=\dd - i \sum_{k=1}^d j_k(\omega_k + \brr{\omega_k})\wedge \cdot$ if $P=\LLl$ and $\delta_j := \dbar - i\sum_{k=1}^d j_k\omega_k \wedge \cdot$ if $P=\dbar_b$ (we omit the degree $q$ in the notation for simplicity). Since $P_j$ form elliptic complexes for every fixed $j\in \Z^d$, we conclude that $\widehat{\alpha}(j) \in C^\infty(M;\Lambda^q \C T^\ast M)$ for every $j\in \Z^d$.

\vspp

We now turn to main result of this section.
\begin{Thm}\label{smooth_criterion} Let $P$ denote either $\LLl$ or $\dbar_b$. Let $\alpha \in X_P^q$, for $q\geq 1$. Then, the following are equivalent:
\begin{enumerate}
\item $\alpha \in C^\infty(M\times \T^d;\Lambda^q\C T^\ast M)$.
\item For every $A\in \Z_+$,
\[
\sup_{j\in \Z^d} (1+|j|)^A \left\|\widehat{\alpha}(j)\right\|_{L^\infty(M)}<\infty.
\]
\item For every $A\in \Z_+$,
\[
\sup_{j\in \Z^d}(1+|j|)^A \left\|\widehat{\alpha}(j)\right\|_{L^2(M)} < \infty.
\]
\end{enumerate}

\end{Thm}
\begin{proof} We first adress the case $P=\LLl$. It is clear that $1) \implies 2) \implies 3)$. To show $3) \implies 1)$, it is enough to show that for every $s,A\in \Z_+$,
\begin{equation}\label{induct_est}
\sup_{j\in \Z^d}(1+|j|)^A \|\widehat{\alpha}(j)\|_{W^s(M)}<\infty,
\end{equation}
where $W^{s}(M)$ is the usual $L^2$ Sobolev space of order $s$. We shall prove \eqref{induct_est} by induction on $s\in \Z_+$.  If $s=0$, this is just our hypothesis $3)$. Assume now that this estimate holds for $s\in \Z_+$. Let $f:=\LLl \alpha \in C^\infty(M\times \T^d;\Lambda^{q+1}\C T^{\ast}M)$ and $f^{\ast}:=\LLl^{\ast}\alpha \in C^\infty(M\times \T^d;\Lambda^{q-1}\C T^{\ast}M)$. We can write $f = \sum_{j\in \Z_d}f_j e^{i\langle j,\theta \rangle}$ and $f^{\ast}=\sum_{j\in \Z^d}f^{\ast}_j e^{i\langle j,\theta \rangle}$, where $f_j = \LLl_j \widehat{\alpha}(j)$ and $f^{\ast}_j = \LLl_j^{\ast}\widehat{\alpha}(j)$ for all $j\in \Z^d$. Then, $\widehat{\alpha}(j)$ satisfies the two following equations:
\[
\begin{cases}
\dd \widehat{\alpha}(j) = i\sum_{k=1}^dj_k (\omega_k+\brr{\omega_k})\wedge \widehat{\alpha}(j)+ f_j \\
\dd^{\ast} \widehat{\alpha}(j) = i\sum_{k=1}^d j_k(\omega_k + \brr{\omega_k})\iprod \widehat{\alpha}(j) + f^{\ast}_j.
\end{cases}
\]
Since the de Rham complex is elliptic, there is\footnote{Indeed, since the Laplacian $\Delta = \dd^{\ast}\dd + \dd \dd^{\ast}$ is an elliptic operator of order $2$ (acting on $q$-forms), we have the fundamental estimate $\|\beta\|_{s+1}\leq C(\|\dd^{\ast}\dd \beta + \dd \dd^{\ast}\beta\|_{s-1}+\|\beta\|_{0})$, which immediately implies the result using the continuity of $\dd$ and $\dd^{\ast}$ from $W^{s}$ to $W^{s-1}$.} a constant $C>0$ such that 
\begin{equation}\label{elliptic_ineq}
\|\beta\|_{W^{s+1}(M)} \leq C \left(\|\dd \beta\|_{W^{s}(M)}+\|\dd^{\ast}\beta\|_{W^s(M)}+\|\beta\|_{L^2(M)} \right),\,\,\beta \in C^\infty(M;\Lambda^{q}\C T^{\ast}M).
\end{equation}
Applying it with $\beta=\widehat{\alpha}(j)$ yields
\begin{align*}
\|\widehat{\alpha}(j)\|_{W^{s+1}(M)} \leq C_1\left((1+|j|)\|\widehat{\alpha}(j)\|_{W^s(M)}+\|f_j\|_{W^s(M)}+\|f_j^{\ast}\|_{W^{s}(M)}+\|\widehat{\alpha}(j)\|_{L^2(M)} \right),
\end{align*}
where $C_1>0$ is independent of $j$. Since $\LLl \alpha$ and $\LLl^{\ast} \alpha$ are both smooth, this implies \ref{induct_est} for $s+1$. The case where $P=\dbar$ is done in a similar way (observing that the Dolbeault complex is also elliptic on $M$).
\end{proof}

\section{Global hypoellipticity}

\n {\bf (A)} We shall now apply the techniques developed in the previous section to relate regularity properties of $\LLl$ and $\dbar_b$. The main property that will concern us is the following:
\begin{Def} Let $(N,g)$ be a Riemannian manifold and $E_1,E_2,E_3$ be smooth vector bundles over $N$, endowed with euclidean (or hermitian) metrics. Let $\mathcal{C}$ be a complex of linear partial differential operators
\[
\mathcal{C}: \Dli(N;E_1) \xrightarrow{Q} \Dli(N;E_2) \xrightarrow{P}\Dli(N;E_3)
\]
with smooth coeficients. We say $\mathcal{C}$ is \textit{globally hypoelliptic} if every distributional section $u\in \Dli(N;E_2)$ such that $Pu \in C^\infty(N;E_3)$ and $Q^{\ast}u \in C^\infty(N;E_1)$ belongs to $C^\infty(M;E_1)$ (here, $Q^\ast$ denotes the formal adjoint of $Q$ with respect to the metric structures).
\end{Def}
In our setting, we use the following terminology:
\begin{Def} Let $\dbar_b$ (respectively, $\LLl$) be the differential complex defined by \ref{cr_complex} (respectively, \ref{real_complex}), and fix a metric on $M\times \T^d$ of the form $g_{M}\oplus g_{\T^d}$, where $g_M$ is a hermitian (respectively, Riemannian) metric on $M$ and $g_{\T^d}$ is the flat metric. We say this complex is \textit{globally hypoelliptic} in degree $q$ ($1\leq q \leq 2n$) if the complex 
\[
\Dli\left(M\times \T^d;\Lambda^{q-1}\C T^{\ast}M\right)\xrightarrow{(\dbar_b)_{q-1}}\Dli\left(M\times \T^d;\Lambda^{q}\C T^{\ast}M\right) \xrightarrow{(\dbar_b)_{q}} \Dli\left(M\times \T^d;\Lambda^{q+1}\C T^{\ast}M\right)
\]
is globally hypoelliptic (respectively, if the complex given by $\{\LLl_{q-1},\LLl_q\}$ is globally hypoelliptic). We say $\dbar_b$ (respectively, $\LLl$) is globally hypoelliptic in degree $0$ if it is globally hypoelliptic in the usual sense (see, for example, \cite{BCM}).
\end{Def}
It is not clear to what extent this definition depends on the Riemannian metric on the manifold $M\times \T^d$. We can show, however, the following characterization:
\begin{Thm}\label{main_thm_1} Let $u\in \Dli(M\times \T^d;\Lambda^q \C T^{\ast}M)$. Let $P$ stand for either $\LLl$ or $\dbar_b$. Then, the following are equivalent:
\begin{enumerate}
\item $P_{q} u \in C^\infty(M\times \T^d;\Lambda^{q+1} \C T^{\ast}M)$ and $P_{q-1}^{\ast}u \in C^\infty(M\times \T^d;\Lambda^{q-1}\C T^{\ast}M)$.
\item $\Box_{P}^{q}u \in C^\infty(M\times \T^d;\Lambda^{q}\C T^{\ast}M)$.
\end{enumerate}
In particular, the complex $\LLl$ (respectively, $\dbar_b$) is globally hypoelliptic in degree $q$ if and only if the corresponding Laplacian
\[
\Box_{\LLl}^q:=\LLl_{q}^{\ast}\LLl_q + \LLl_{q-1}\LLl_{q-1}^{\ast}: \Dli\left(M\times \T^{d};\Lambda^{q}\C T^{\ast}M \right)\to \Dli\left(M\times \T^{d};\Lambda^{q}\C T^{\ast}M \right)
\]
is globally hypoelliptic in the usual sense (respectively, the Laplacian $\Box_{\dbar_b}^q$ is globally hypoelliptic). If $q=0$, the Laplacians are defined by $\Box_{\LLl}^0 := \LLl_0^{\ast}\LLl_0$ (respectively, $\Box_{\dbar_b}^{0} := (\dbar_b)_0^{\ast}(\dbar_b)_0)$.
\end{Thm}
\begin{proof} Fix a degree $0\leq q \leq 2n$ and consider the operator $\dbar_b$ (the proof for $\LLl$ is the same). It is clear that if $\Box_{\dbar_b}^q$ is globally hypoelliptic, then $\dbar_b$ is globally hypoelliptic in degree $q$. To show the converse, assume $\dbar_b$ is globally hypoelliptic in degree $q$ and let $u\in \Dli(M\times \T^d;\Lambda^{q}\C T^{\ast}M)$ be a distributional sectional such that $\Box_{\dbar_b}^q u \in C^\infty(M\times \T^{d};\Lambda^{q} \C T^{\ast}M)$. Write the Fourier decomposition
\[
u=\sum_{j\in \Z^{d}}\widehat{u}(j)e^{i\langle j,\theta \rangle},\,\,\widehat{u}(j) \in C^\infty(M;\Lambda^{q}\C T^{\ast}M).
\]
(observe that the Fourier coefficients of $u$ are smooth since $(\Box_{\dbar_b}^{q})_j \widehat{u}(j) \in C^\infty(M;\Lambda^{q}\C T^{\ast}M)$ for every $j \in \Z^d$, and these are elliptic complexes for every fixed $j \in \Z^d$). We shall first prove that $(\dbar_b)_{q}^{\ast}(\dbar_b)_q u \in C^\infty(M\times \T^{d};\Lambda^q \C T^{\ast}M)$. Since
\[
(\dbar_b)_{q-1}^{\ast}(\dbar_b)_{q}^{\ast}(\dbar_b)_q u = 0 \text{ and }(\dbar_b)_q(\dbar_b)_{q}^{\ast}(\dbar_b)_q = (\dbar_b)_q \Box_{\dbar_b}^{q}u \in C^\infty(M\times \T^d;\Lambda^{q}\C T^{\ast}M),
\]
we have $(\dbar_b)_{q}^{\ast}(\dbar_b)_q u \in X^{q}_{\dbar_b}$. By Theorem \ref{smooth_criterion}, we have to estimate the $L^{2}$ norm of the Fourier coefficients of $(\dbar_b)_{q}^{\ast}(\dbar_b)_q u $, which are given by $\delta_j^{\ast}\delta_j \widehat{u}(j)$ (we omit the degree $q$ from the notation for simplicity). Fix $j\in \Z^d$. Then,
\[
\left\|\delta_j^{\ast}\delta_j \widehat{u}(j)\right\|_{L^2(M)}^2 = \llangle \delta_j \widehat{u}(j),\delta_j \left(\Box_{\delta_j}\widehat{u}(j)\right)\rrangle.
\]
Since $M\times \T^d$ is a compact manifold, there is a number $s\in \R$ such that $\dbar_b^q u \in H^{s}(M\times \T^d;\Lambda^{q}\C T^{\ast}M)$. In particular, we have
\[
\sum_{j\in \Z^d} (1+|j|)^{s}\|\delta_j \widehat{u}(j)\|_{W^s(M)}<\infty.
\]
We conclude then, from the generalized Cauchy-Schwarz inequality (see, for example, Proposition A.1.1 in \cite{Folland-Kohn})
\begin{align*}
\left\|\delta_j^{\ast}\delta_j \widehat{u}(j)\right\|_{L^2(M)}^2 \leq \|\delta_j \widehat{u}(j)\|_{W^{s}(M)} \|\delta_j (\Box_{\delta_j}\widehat{u}(j))\|_{W^{-s}(M)}.
\end{align*}
Therefore, given $A\in \Z_+$,
\begin{align*}
(1+|j|)^{2A}\left\|\delta_j^{\ast}\delta_j \widehat{u}(j)\right\|_{L^2(M)}^2  \leq \left\{(1+|j|)^s \|\delta_j \widehat{u}(j)\|_{W^s(M)} \right\} \left\{(1+|j|)^{2A-s} \|\delta_j (\Box_{\delta_j}\widehat{u}(j))\|_{W^{-s}(M)}\right\}
\end{align*}
Since $\dbar_b (\Box_{\dbar_b}^q u) $ is smooth, the supremum of the expression above is finite over $j\in \Z^d$. We conclude by theorem \ref{smooth_criterion} that $(\dbar_b)^{\ast}_q (\dbar_b)_q u$ is smooth. In the same way one proves that $(\dbar_b)_{q-1}(\dbar_b)^{\ast}_{q-1}u$ is also smooth. We obtain then that $(\dbar_b)_q u \in X^{q+1}_{\dbar_b}$ (also $(\dbar_b)^{\ast}_{q-1}u \in X^{q-1}_{\dbar_b}$). Now, writing
\[
\|\delta_j \widehat{u}(j)\|_{L^2(M)}^2 = \llangle \widehat{u}(j),\delta_j^{\ast}\delta_j \widehat{u}(j)\rrangle
\]
and using smoothness of $(\dbar_b)^{\ast}_q (\dbar_b)_q u$, we obtain again by Theorem \ref{smooth_criterion} that $(\dbar_b)_q u$ is smooth. In the same way, $(\dbar_b)^{\ast}_{q-1}u$ is smooth. Since $\dbar_b$ is globally hypoelliptic in degree $q$, the section $u$ is smooth. The result is proved.

\end{proof}
We can now state the main result of this section:
\begin{Cor}\label{73} Let $0\leq q \leq 2n$ and $u\in \Dli(M\times \T^d;\Lambda^{q}\C T^{\ast}M)$. Assume that $M$ is a balanced manifold in degree $q$ (see definition \ref{def1}). Then, the following are equivalent:
\begin{enumerate}
\item $(\dbar_b)_q u \in C^\infty(M\times \T^d;\Lambda^{q+1}\C T^\ast M)$ and $(\dbar_b)^{\ast}_{q-1}u \in C^\infty(M\times \T^d;\Lambda^{q-1}\C T^{\ast}M)$.
\item $\LLl_q u \in C^\infty(M \times \T^d;\Lambda^{q+1}\C T^\ast M)$ and $\LLl_{q-1}^{\ast}u \in C^\infty(M\times \T^{d};\Lambda^{q-1}\C T^{\ast}M)$. 
\end{enumerate}
In particular, $\dbar_b$ is globally hypoelliptic in degree $q$ if and only if $\LLl$ is globally hypoelliptic in degree $q$.
\end{Cor}
In many interesting cases, a characterization for global hypoellipticity (in degree $0$) is known for $\LLl$. For instance, using Theorem 7.3 in \cite{Araujo2022} we obtain the following
\begin{Cor}\label{Mahler} Assume that $M$ is balanced in degree $0$. Then, $\dbar_b$ is globally hypoelliptic in degree $0$ if and only if the system $\bm{\omega+\overline{\omega}}=(\omega_1+\overline{\omega_1},\ldots,\omega_d+\overline{\omega_d})
$ is neither rational nor Liouville.
\end{Cor}
\begin{Rem} We refer to the work \cite{Araujo2022} for the definition of rational and Liouville systems. We also remark that this result (concerning global hypoellipticity of $\LLl$ in degree $0$) for $d=1$ was obtained by \cite{BCM}.
\end{Rem}

\section{Global solvability}

\n {\bf(A)} We shall now compare the cohomology of the complexes \eqref{real_complex} and \eqref{cr_complex}. We use the following notation: 
\[
\mathcal{H}^{q}(M\times \T^d;P) = \left\{u\in C^\infty(M\times \T^d;\Lambda^q \C T^{\ast}M);\,\Box_{P}u=0 \right\},
\]
where $P\in\{\LLl,\dbar_b\}$. We have the natural (injective) map
\[
i^{P}_q:\mathcal{H}^{q}(M\times \T^d;P) \to H^q(M\times \T^d;P),
\]
sending a $P$-harmonic smooth form to its cohomology class (again, for $P\in \{\LLl,\dbar_b\}$). We have the following result:
\begin{Prop} Let $P\in \{\LLl,\dbar_b\}$ and $0\leq q \leq 2n$. Then, the following are equivalent:
\begin{enumerate}
\item $H^{q}(M\times \T^d;P)$ is a Fréchet space.
\item $i^{P}_q$ is a topological isomorphism.
\end{enumerate}
\end{Prop}
\begin{proof} The implication $2)\implies 1)$ is immediate. For the converse direction, it remains to see $i^P_q$ is surjective. Let $[u]\in H^q(M\times \T^d;P)$ be a cohomology class, and choose a representative $u\in C^\infty(M\times \T^d;\Lambda^q\C T^{\ast}M)$ such that $Pu=0$. We write the Fourier series decomposition
\[
u=\sum_{j\in \Z^d}u_j e^{ij\theta},
\]
where $u_j\in C^\infty(M;\Lambda^q \C T^{\ast}M)$ is such that $P_j u_j = 0$ for all $j\in \Z^d$ ($P_j$ is either $D_j$ or $\delta_j$, as in section $4$). Recalling these complexes are elliptic, there is a unique $P_j$-harmonic representative $v_j \in C^\infty(M;\Lambda^q \C T^{\ast}M)$ in the $P_j$-cohomology class of $u_j$, i.e., $u_j-v_j = P_j w_j$ for some $(q-1)$-form $w_j \in C^\infty(M;\Lambda^{q-1} \C T^{\ast}M)$ and $\Box_{P_j}v_j = 0$. We also know that $\|v_j\|_{L^2(M)} \leq \|u_j\|_{L^2(M)}$ for all $j\in \Z^d$. Therefore, by Theorem \ref{smooth_criterion}, we conclude that 
\[
v=\sum_{j\in \Z^d}v_j e^{ij\theta} 
\]
defines an element in $\mathcal{H}^q(M\times \T^d;P)$. It remains to see that $i^{P}_q(v)=[u]$. Indeed, for every $N\geq 1$, we have
\[
\sum_{|j|\leq N}(u_j-v_j)e^{ij\theta}=P \left(\sum_{|j|\leq N}w_j e^{ij\theta} \right).
\]
This sequence of elements in the (closed) range of $P$ converges to $u-v$, which must be also in the range of $P$ by $(1)$. The result is proved.
\end{proof}
 \begin{Rem} These results do not depend on any metric properties of $M$, only on the ellipticity of the de Rham and Dolbeault's complexes.
\end{Rem}
We move to the main result of this section.
\begin{Thm}\label{comp_cohomol_1} Let $M$ be balanced in degree $0$. Then, the following are equivalent:
\begin{enumerate}
\item $H^{1}(M\times \T^d;\LLl)$ is a Fréchet space.
\item $H^{1}(M\times \T^d;\dbar_b)$ is a Fréchet space.
\end{enumerate}
Under these conditions, we have $H^{1}(M\times \T^d;\LLl) \simeq H^{1}(M\times \T^d;\dbar_b)$. We also have $H^0(M\times \T^d;\LLl)=H^0(M\times \T^d;\dbar_b)$.
\end{Thm}

The main technical ingredient is the following:
\begin{Prop} Let $M$ be balanced in degree $0$. Then, for every $k\in \Z_+$, there is a constant $C_k>0$ such that
\begin{equation}\label{comp_estimate}
\left\|\dbar_b u\right\|_{W^k} \leq \left\|\LLl u\right\|_{W^k} \leq C_k \left\|\dbar_b u\right\|_{W^k},\,\,\,u\in C^\infty(M\times \T^d).
\end{equation}

\end{Prop}
\begin{proof} The inequalities \ref{comp_estimate} follows from \ref{teo6} for $k=0$, so we assume $k\geq 1$. The first inequality in \ref{comp_estimate} is immediate, since $\dbar_b u$ is the projection of $\LLl u$ onto the $(0,1)$-forms.  To prove the second estimate, we argue by contradiction: assume the estimate is false. Then, we can find $k\geq 1$ and a sequence $u_n \in C^\infty(M\times \T^d)$ such that $\left\|\LLl u_n\right\|_{W^k}=1$ for all $n$ and $\dbar_b u_n \to 0$ in $W^k$. 

By Rellich's lemma, we can (after passing to a subsequence) assume that $\LLl u_n \to v$ in $W^{k-1}$, where $v\in W^{k-1}(M\times \T^d;\C T^{\ast}M)$. Applying $\LLl^{\ast}$ yields (using \ref{teo6})
\[
2\dbar_b^{\ast}\dbar_b u_n = \LLl^{\ast}\LLl u_n \to \LLl^{\ast}v \text{ in }\Dli.
\]
We also have $\dbar_b^{\ast}\dbar_b u_n \to 0$, so by uniqueness of the distributional limit, we conclude that $\LLl^{\ast}v=0$, i.e., $v\in \ker \LLl^{\ast}\cap W^{k-1}(M\times \T^d;\C T^{\ast}M)$. Since $k\geq 1$, we have $v\in L^{2}$ and since the range of $\LLl$ is orthogonal to the kernel of $\LLl^{\ast}$, we conclude that $v=0$, which contradicts the fact that $\|\LLl u_n\|_{W^k}=1$ for all $n$. 
\end{proof}
A consequence of this proposition is the following
\begin{Cor} Assume $M$ is balanced in degree $0$. Then, the following are equivalent:
\begin{enumerate}
\item $\LLl: C^\infty(M\times \T^d) \to C^\infty(M\times \T^d;\C T^{\ast}M)$ has closed range.
\item $\dbar_b: C^\infty(M\times \T^d) \to C^\infty(M\times \T^d;\C T^{\ast}M)$ has closed range.
\end{enumerate}

\end{Cor}
\begin{proof} We apply the following characterization from \cite{Kothe}, page 18: a continuous linear map $T:E\to F$ between Fréchet spaces has closed range if, and only if, the following holds:
\begin{align}\label{closed_range_charact}
&\text{For every sequence }(u_n)\text{ in }E\text{ such that }Tu_n\to 0,\\
&\text{ there is a sequence }v_n\in E\text{ such that }Tu_n=Tv_n\text{ and }v_n\to 0. \nonumber
\end{align}
Assume $\LLl$ has closed range and let $u_n \in C^\infty(M\times \T^d)$ be such that $\dbar_b u_n \to 0$. From \ref{comp_estimate}, we have $\LLl u_n \to 0$ in the $C^\infty$ topology. From \ref{closed_range_charact}, there is a sequence $v_n \in C^\infty(M\times \T^d)$ such that $\LLl u_n=\LLl v_n$ and $v_n \to 0$. However, the kernel of $\LLl$ and $\dbar_b$ are equal in $C^\infty(M\times \T^d)$ (again from \ref{comp_estimate}), so \ref{closed_range_charact} implies that $\dbar_b$ has closed range. The argument in the other direction is identical.
\end{proof}
Now, recalling that $H^1(M\times \T^d;\LLl)$ is Fréchet if and only if $\LLl$ has closed range (same for $\dbar_b$), the statement \ref{comp_cohomol_1} follows immediately. 

\vspp

In recent work \cite{Araujo2023}, a characterization for closedness of the range of $\LLl$ was obtained in terms of a diophantine condition on the forms. Using this result, we can state the following

\begin{Cor} Assume $M$ is balanced in degree $0$. Then, $\dbar_b: C^\infty(M\times \T^d) \to C^\infty(M\times \T^d;\C T^{\ast}M)$ has closed range if and only if the collection $\bm{\omega+\brr{\omega}}=(\omega_1+\brr{\omega_1},\ldots,\omega_d + \brr{\omega_d})$ is weakly non-simultaneously approximable. 
\end{Cor}
\begin{Rem} We refer to the paper \cite{Araujo2023} for the definition of a weakly non-simultaneously approximable system.
\end{Rem}

We say $\LLl$ (respectively, $\dbar_b$) is \textit{globally solvable} if $H^1(M\times \T^d;\LLl)=0$ (respectively, $H^{1}(M\times \T^d;\dbar_b)=0$). In view of the previous results, we can state the following
\begin{Cor} Let $M$ be balanced in degree $0$. Then, the following are equivalent:
\begin{enumerate}
\item $\LLl$ is globally solvable.
\item $\dbar_b$ is globally solvable.
\end{enumerate}
\end{Cor}

\section*{Acknowledgments}
\n  The first author was partially supported by Conselho Nacional de Desenvolvimento Cient\'ifico e Tecnol\'ogico-CNPq, Grant 303945/2018-4 and S\~ao Paulo Research Foundation (FAPESP), Grant 2018/14316-3. The second author was partially supported by Conselho Nacional de Desenvolvimento Científico e Tecnológico-CNPq, Grant 400739/2022-4 and by the Czech Science Foundation, project GC22-15012J.

\printbibliography

\vskip 0.3in

\begin{minipage}[b]{8 cm}
	{\bf Paulo D. Cordaro} ({\it \small Corresponding author})
	\\ 
	University of S\~ao Paulo\\ 
	S\~ao Paulo, Brazil\\
	E-mail: {\it cordaro@ime.usp.br}
\end{minipage}

\vskip 0.3in
\begin{minipage}[b]{8 cm}
	{\bf Vinícius Novelli}\\
	Masaryk University\\
        Brno, Czechia\\
	E-mail: {\it novelli@math.muni.cz}
\end{minipage}

\end{document}